\documentclass[12pt]{amsart}
\usepackage[scale=.86]{geometry}
\usepackage{geometry} 
\usepackage{graphicx}
\usepackage{bbm}
\usepackage{mathrsfs}
\theoremstyle{plain}

\usepackage{CJK}
\usepackage{amsfonts,amssymb,comment,amsmath,mathrsfs,multirow,booktabs}

 \usepackage{hyperref}  

\hypersetup{
	colorlinks = true,  
	linkcolor = blue,  
	citecolor = blue,   
	urlcolor  = red     
}

\usepackage{color,latexsym,amsfonts}
 \newtheorem{theorem}{Theorem}[section]
 \newtheorem{lemma}[theorem]{Lemma}
 \newtheorem{corollary}[theorem]{Corollary}
 \newtheorem{proposition}[theorem]{Proposition}
 \newtheorem{example}[theorem]{Example}
 \newtheorem{Definition}[theorem]{Definition}
 \newtheorem{remark}[theorem]{Remark}
 \newtheorem{condition}[theorem]{Condition}

 \def\beqlb{\begin{eqnarray}}\def\eeqlb{\end{eqnarray}}
 \def\beqnn{\begin{eqnarray*}}\def\eeqnn{\end{eqnarray*}}

 \def\mcr{\mathscr}\def\mbb{\mathbb}\def\mbf{\mathbf}

 \def\<{\langle}\def\>{\rangle}

 \def\ar{&\!\!}\def\crr{\\}

 \def\eqref#1{{\rm(\ref{#1})}}

 \def\proof{\noindent{\it
 Proof.~}}\def\qed{\hfill$\Box$\medskip}

\def\e{{\mbox{\rm e}}}
\def\<{\left<}\def\>{\right>}

 \def\mcr{\mathscr}
 \def\mbb{\mathbb}
  \def\mbf{\mathbf}
\newcommand{\dd}{\mathrm{d}}

\newfam\msbmfam\font\tenmsbm=msbm10\textfont
\msbmfam=\tenmsbm\font\sevenmsbm=msbm7
\scriptfont\msbmfam=\sevenmsbm

\def\<{\left<}\def\>{\right>}

\def\({\left(}\def\){\right)}

\usepackage[dvipsnames]{xcolor}

\begin{document}
	
	
	
	\centerline{\large\textbf{Boundary behavior at infinity for  }}
	
	\smallskip
	
	\centerline{\large\textbf{simple exchangeable fragmentation-coagulation process}}

	\smallskip
	\centerline{\large\textbf{ in critical slow regime}}
	
	\bigskip
	
	\centerline{Lina Ji\,\footnote{Supported by NSFC grant (No. 12301167), Guangdong Basic and Applied Basic Research Foundation (No. 2022A1515110986) and Shenzhen National Science Foundation (No. 20231128093607001).} and Xiaowen Zhou\,\footnote{Supported by Natural Sciences and Engineering Research Council of Canada (RGPIN-2021-04100).}}
	
	\medskip
	
	\centerline{\it  Joint Research Center of Applied Mathematics, }

    \centerline{\it Shenzhen MSU-BIT University, Shenzhen 518172, P.R. China.}

    \centerline{\it Department of Mathematics and Statistics, Concordia University}

    \centerline{\it 1455 De Maisonneuve Blvd. W., Montreal, Canada.}
	
	\centerline{E-mails: \tt jiln@smbu.edu.cn, xiaowen.zhou@concordia.ca}
	
	\bigskip
	
	{\narrower{\narrower

\centerline{\textbf{Abstract}}

\bigskip
For a critical simple exchangeable fragmentation-coagulation process in slow regime where the coagulation rate and fragmentation rate are of the same order, we show that there exist phase transitions for its boundary behavior at infinity depending on the asymptotics of the difference between the two rates, and find rather sharp conditions for different boundary behaviors.

\bigskip

\medskip


\par}\par}

\bigskip
	
	\section{Introduction}

Fragmentation and coagulation  process appears in various physical and biological models. Intuitively, it describes a particle system in which  particles can merge to form larger clusters and can also break into smaller ones. The exchangeable
fragmentation-coagulation processes (EFC-processes for short) were introduced by Berestycki \cite{B04} as partition-valued processes. Roughly speaking, an EFC process, defined on the space of partitions, evolves in continuous time and combines the dynamics of coagulations and homogeneous fragmentations.  
We refer to Bertoin \cite{B06}, Pitman \cite{P06} and references therein  for an introduction to exchangeable fragmentations and coalescents.

For $n \in \mbb{N}_+ := \{1, 2, \cdots\}$, let $\mathcal{P}_n$ denote the collection of partitions of $[n]:= \{1, \cdots, n\}$, where a partition $\pi\equiv (\pi_i)\in\mathcal{P}_n$ consists of disjoint subsets $\pi_i$s of $[n]$ ordered by their least elements and satisfying $\cup_i\pi_i=[n]$. Let $\mathcal{P}_\infty$ be the set of partitions on $\mbb{N}_+.$
For any $\pi \in \mathcal{P}_\infty$ with $\pi = (\pi_1, \pi_2, \cdots)$, let $\pi_{|[n]}$ be the restricted partition $(\pi_i\cap[n], i \ge 1).$ Then $\pi_{|[n]} \in \mathcal{P}_n$. A $\mathcal{P}_\infty$-valued Markov process $(\Pi(t), t \ge 0)$ is an {\it EFC process} if satisfies:
\begin{itemize}
	\item It is exchangable, i.e., for any time $t \ge 0$, the random partition $\Pi(t)$ of $\mbb{N}_+$ has a law invariant under the permutations on $[n]$ for any $n \in \mbb{N}_+$;
	\item For any $n \in \mbb{N}_+,$ the restriction $(\Pi(t)_{|[n]}, t \ge 0)$ is a c\`adl\`ag Markov chain, taking values on $\mathcal{P}_n$, which can only evolve by fragmentation of one block or by coagulation.
\end{itemize}
See \cite[Definition 1]{B04}. For definitions of fragmentation and coagulation, we refer readers to \cite[Definitions 3.1 and 4.2]{B06}.
It is worth noting that any block in an exchangeable random partition of $\mbb{N}_+$ is either a singleton or an infinite block, see \cite[Subsection 4.2]{B04} and Foucart \cite[Subsection 2.1]{F22}. An EFC process is called {\it simple} if it excludes simultaneous multiple collisions and instantaneous coagulation of all blocks, and its fragmentation event occurs at a finite total rate and exclusively involves singleton-free partitions; see
\cite[Definition 2.9]{F22}. In such processes, fragmentation occurs at a finite rate without generating singletons. Consequently, the process remains {\it proper} at all times, i.e., an exchangeable partition contains no singleton block in the sense of \cite[Chapter 2.3]{B06}.

The simple EFC-proceses can be seen as a generalization of $\Lambda$-coalescents defined by Pitman \cite{P99} and Sagitov \cite{S99}, see also Berestycki \cite{B09}.
The $\Lambda$-coalescent is an exchangeable coalescent with multiple collisions.
For any $2 \le k \le n$, the block counting process of $\Lambda$-coalescent jumps from $n$ to $n-k + 1$ at rate  
\beqlb\label{lambda}
\lambda_{n,k}\ar:=\ar  \int_{[0,1]} x^{k-2}(1 - x)^{n-k}\Lambda(\dd x),
\eeqlb
where $\Lambda$ is a finite measure on $[0, 1]$. Its coming down from infinity property has been extensively studied in the 2000s, see Berestycki et al. \cite{BBL10}, Schweinsberg \cite{S00} and Limic and Talarczyk \cite{LT15}. In particular,
a necessary and sufficient condition for a $\Lambda$-coalescent to come down from
infinity was given by \cite{S00}: define
for any $n \ge 2$,
\beqlb\label{Phi}
\Phi_\Lambda(n) := \sum_{k = 2}^n  \binom{n}{k} \lambda_{n,k}(k - 1),
\eeqlb
the $\Lambda$-coalescent comes down from infinity if and only if $\sum_{n=2}^\infty \frac{1}{ \Phi_\Lambda(n)} < \infty.$

In a parallel setting without coalescence, the block counting process with finite initial value jumps from $n$ to $n + k$ at rate $n\mu(k)$ for any $k \in \bar{\mbb{N}}_+  := \mbb{N}_+\cup\{\infty\}$ and a finite measure $\mu$ on $\bar{\mbb{N}}_+$ satisfying $\mu(\infty) = 0$, which guarantees that fragmentation cannot split a block into infinitely many sub-blocks. We shall call $\mu$ the {\it splitting measure}. Therefore, the block counting process can also be interpreted as a continuous-state Markov branching process with non-decreasing paths, which can diverge to infinity in finite time (see Athreya and Ney \cite{AN72}) causing explosion.

When both fragmentation and coalescence are considered, the sample paths of the block-counting process are no longer monotonic, and new phenomena may arise due to the interplay between the fragmentation and the coalescence.

Let $\Pi:= (\Pi(t), t \ge 0)$ denote a simple EFC process and  $N:= (N_t, t \ge 0)$ denote the corresponding block-counting process, where  $N_t := \sharp \Pi(t)$ is the total  number of non-empty blocks in $\Pi(t)$. The boundary behavior of the EFC process are closely related to that of the corresponding block-counting process $N$. Now we present some definitions concerning the boundary behavior of $N$.

\begin{Definition}		
	Let $(N_t, t \ge 0)$ be the block counting process with $N_0 =\infty$. It comes down from infinity if $\mbf{P}\{N_t < \infty\ \text{for\ some}\ t > 0\} = 1.$ Otherwise, it stays infinite.
\end{Definition}
\begin{Definition}
	Let $(N_t, t \ge 0)$ be the block counting process with $N_0 < \infty$. It explodes if $\mbf{P}\{N_t = \infty\ \text{for\ some}\ t > 0\} = 1$. Otherwise, it does not explode.
\end{Definition}	
\begin{Definition}
	Let $(N_t, t \ge 0)$ be the block counting process.
	\begin{itemize}
		\item $\infty$ is an entrance boundary for the process $N$ if $N$ does not explode and comes down from infinity;	
		\item $\infty$ is an exit boundary for the process $N$ if $N$ explodes and stays infinite.
	\end{itemize}
\end{Definition}

For the EFC process with binary coalescence, under certain assumptions on fragmentation, the block-counting process with a finite initial value exhibits the same dynamic as a discrete logistic branching process (see \cite[Section 5]{B04} and Lambert \cite[Section 2.3]{L05}). This connection allows for the derivation of a sufficient condition for the process to come down from infinity, as established in \cite[Proposition 15]{B04}. Later, Kyprianou et al. \cite{KPRS17} investigated the boundary behavior at infinity for a ``fast" EFC process with binary coalescence, where fragmentation dislocates each individual block into its constituent singletons at a constant rate.
Recently, the simple EFC process with $\Lambda$-coalescences has been studied in \cite{F22} and Foucart and Zhou \cite{FZ22}, where fragmentation dislocates at finite rate an individual block into sub-blocks of infinite size. A phase transition between a regime in which $N$ comes down from infinity and one in which it stays infinite is established in \cite[Theorem 1.1]{F22}. The conditions for explosion and non-explosion of the process are provided in \cite[Theorems 3.1 and 3.3]{FZ22}. These results are combined to study the nature of the boundary at $\infty$ for slower-varying coalescence and fragmentation mechanisms, which is summarized below.  Note that we have corrected a typo in \cite[Theorem 3.9]{FZ22} on the assumption of $\beta$.

\begin{theorem}[Theorem 3.9 of \cite{FZ22}]\label{FZ3.9}
Assume that $\mu(n) \underset{n \rightarrow \infty}{\sim} b(\log n)^\alpha n^{-2}$ and $\Phi_\Lambda(n) \underset{n \rightarrow\infty}{\sim} dn(\log n)^\beta$ for $b, d, \alpha> 0$ and $\beta > 1$.
\begin{itemize}
\item If $\beta < 1 + \alpha$, then $\infty$ is an exit boundary;
\item If $\beta > 1 + \alpha$, then $\infty$ is an entrance boundary;
\item If $\beta = 1 + \alpha$ and further,
\begin{itemize}
	\item if $d < \frac{b}{1 + \alpha}$, then $\infty$ is an exit boundary;
	\item if $d > \frac{b}{1 + \alpha}$, then $\infty$ is an entrance boundary.
\end{itemize}
\end{itemize}
\end{theorem}

However, the nature of boundary $\infty$ was left open in some critical regimes, such as the case $d = \frac{b}{1 + \alpha}$ with $\beta = 1 + \alpha$ in Theorem \ref{FZ3.9}. In this paper we  further investigate this issue and identify phase transitions for the boundary behaviors when additional conditions are imposed on asymptotics of the differences between $\Phi_\Lambda(n)$ and   
\beqnn
\Phi_\mu(n) := n\sum_{k = 1}^n\mu(k)k, \quad  n=1,2,\ldots.
\eeqnn

Throughout this paper, we assume that the following condition holds:
\begin{condition}\label{C1}
$\mu(n) \underset{n \rightarrow \infty}{\sim} b(\log n)^\alpha n^{-2}$ as $n\to\infty$ for $b, \alpha > 0.$
\end{condition}

Our approach differs  from those in  \cite{F22,FZ22}. We first refine Chen's method  to develop boundary behavior criteria for the block-counting process on $\bar{\mbb{N}}_+$.
Proofs of the main results then ultimately reduce to identifying appropriate test functions for these criteria to reach the best possible results, where localization conditions in the criteria simplify the construction of test functions and  a coupling for EFC process established in \cite{F22}  helps to complete the proof. It is worth  mentioning that the classification criteria was originally used in  Chen \cite{C04,C86a,C86b} and in Meyn and Tweedie \cite{MT93} for Markov chains, and was later applied to stochastic differential equations associated with continuous-state branching processes in Li et al. \cite{LYZ19}, Ma et al. \cite{MYZ21}, Ren et al. \cite{RXYZ22} and Ma and Zhou \cite{MZ23}.

The rest of this paper is organized as follows. In Section 2, we introduce the model and present our main results. In Section 3, we review  known results of the EFC process and provide  estimates for $\Lambda$ and $\mu$. Section 4 presents criteria for the boundary behaviors of the block-counting process on $\bar{\mbb{N}}_+$. Finally, the proof of the main results is provided in Section 5.

We  conclude this section with some notion to be used.  For any partition $\pi\in \mathcal{P}_\infty$, we denote by $\sharp\pi$, its number of non-empty blocks. By convention, if $\sharp\pi < \infty,$ then we set $\pi_j = \emptyset$ for any $j \ge \sharp\pi + 1.$
Given two functions $f, g: \mbb{R}_+ \mapsto \mbb{R}_+$, write $f = O(g)$ if $\limsup f(x)/g(x) < \infty,$ and $f = o(g)$ if $\limsup f(x)/g(x) = 0$, and $f\sim g$ if $\lim f(x)/g(x) = 1.$ The point at which the limits are taken might vary, depending on the context. For any positive functions $f$ and $g$ well defined on $\mathbb{N}_+$, we write $f = O(g),$ $f = o(g)$ and $f(n)\underset{n\rightarrow \infty}{\sim} g(n)$ if $\limsup_{n \rightarrow \infty}f(n)/g(n)< \infty,$ $\limsup_{n \rightarrow \infty}f(n)/g(n) = 0$ and $\lim_{n \rightarrow \infty} f(n)/g(n) = 1$, respectively.  The $\log$ functions appeared in this paper are all with base $e$. Moreover, the iterated logarithm $\log^{(m)} n \ (m = 0, 1, 2, \cdots)$ is defined as applying the logarithm function recursively $m$ times to $n$, i.e., $\log(\log(\cdots \log(n)))$ with $m$ nested logarithms, where $\log^{(0)} n = n$. For any real number $x$, $\lfloor x\rfloor$ denotes the greatest integer less than or equal to $x$.  We use $C$ to denote a positive constant whose value may change from line to line.

\section{Simple EFC processes and main results}\label{s2}

\subsection{Introduction of simple EFC processes}\label{ss21}
Simple EFC processes are Feller processes with state space $\mathcal{P}_\infty$, which are characterized in law by two $\sigma$-finite exchangeable measures on $\mathcal{P}_\infty$, $\mu_\text{Coag}$ and $\mu_\text{Frag}$, the measures of coagulation and fragmentation, respectively. We briefly recall the Poisson construction of simple EFC processes with given coagulation and fragmentation measures.

Let $(\Omega, \mcr{F}, (\mcr{F}_t)_{t \ge 0}, \mbf{P})$ be a filtered probability space satisfying the usual conditions. On this stochastic basis, we consider two independent Poisson point processes
\beqnn
\text{PPP}_\text{C} = \sum_{t > 0}\delta_{(t, \pi^c)}\quad \text{and}\quad \text{PPP}_\text{F} = \sum_{t > 0}\delta_{(t, \pi^f, k)}
\eeqnn
defined, respectively, on $\mbb{R}_+\times\mathcal{P}_\infty$ and $\mbb{R}_+\times\mathcal{P}_\infty\times\mbb{N}_+$ with intensity $\dd t \otimes \mu_{\text{Coag}}(\dd \pi)$ and $\dd t \otimes\mu_\text{Frag}(\dd \pi)\otimes \sharp(\dd k)$, where $\sharp$ is the counting measure on $\mbb{N}_+$. Here we assume that $\mu_\text{Coag}(\dd \pi)$ is supported on partitions containing more than one block, with exactly one non-singleton block and all other blocks being singletons. Moreover, $\mu_\text{Frag}(\dd \pi)$ is supported on partitions with infinite blocks satisfying $\mu_\text{Frag}(\mathcal{P}_\infty) < \infty$. Let $\Pi(0):= \{\Pi_1(0), \Pi_2(0), \cdots\}$ be an proper and exchangeable random partition independent of $\text{PPP}_\text{C}$ and $\text{PPP}_\text{F}$. For any $n \ge 1$, we set $\Pi^{[n]}(0) = \pi_{|[n]}$. Then the sequence of processes $\{(\Pi^{[n]}(t), t \ge 0): n \ge 1\}$ can be constructed as follows:
\begin{itemize}
\item Coalescence: for each $n \ge 1,$ at an atom $(t, \pi^c)$ of $\text{PPP}_\text{C}$ such that $\pi^c_{|[n]}$ only has one non-singleton block:
\beqnn
\Pi^{[n]}(t) =\text{Coag}(\Pi^{[n]}(t-), \pi^c_{|[n]}),
\eeqnn
where for any partitions $\pi, \pi^c$, $\text{Coag}(\pi,\pi^c) := \{\cup_{j\in\pi_i^c}\pi_j: i \ge 1\}$;
\item Fragmentation: for each $n \ge 1,$ at an atom $(t, \pi^f, k)$ of $\text{PPP}_\text{F}$ such that $\pi^f_{|[n]}$ has at least two non-empty blocks,
\beqnn
\Pi^{[n]}(t) = \text{Frag}(\Pi^{[n]}(t-), \pi^f_{|[n]}, k),
\eeqnn
where for any partitions $\pi, \pi^f$, $ \text{Frag}(\pi,\pi^f,k) := \{\pi_k\cap\pi_i^f, i \ge 1; \pi_\ell, \ell \neq k\}^\downarrow$. Here $\{\cdots\}^\downarrow$ means that   blocks in the partition are ordered by their least elements.
\end{itemize}
The processes $(\Pi^{[n]}(t), t \ge 1)_{n \ge 1}$ are compatible in the sense that for any $m \ge n \ge 1$,
\beqnn
(\Pi^{[m]}(t)_{|[n]}, t \ge 0) = (\Pi^{[n]}(t), t \ge 0).
\eeqnn
This ensures the existence of a process $(\Pi(t), t \ge 0)$ on $\mathcal{P}_\infty$ such that for all $n \ge 1$,
\beqnn
(\Pi(t)_{|[n]}, t \ge 0) = (\Pi^{[n]}(t), t \ge 0).
\eeqnn
The process $\Pi = (\Pi(t), t \ge 0)$ is a simple EFC-process started from $\Pi(0)$, see \cite{F22} and \cite{FZ22}. Moreover, $(\Pi(t), t \ge 0)$ is a general EFC-process by relaxing the restriction on $\mu_\text{Coag}$ and $\mu_\text{Frag}$, see \cite[Subsection 3.2]{B04}.

Note that the merger of multiple blocks into a single block is possible for simple EFC processes. However, under the assumption on $\mu_\text{Coag}$, simultaneous multiple mergers can not occur. Then the $\Lambda$-coalescent can be constructed by this  coagulation measure $\mu_\text{Coag}$, whose probability law  is characterized  by  jump rates of its restrictions, namely, the sequence $(\lambda_{n,k}, 2 \le k \le n)_{n \ge 2}$ defined by  
\beqnn
\lambda_{n,k}\ar:=\ar \mu_\text{Coag}\{\pi;\ \text{the\ non-singleton\ block\ of}\ \pi_{|[n]}\ \text{has}\ k\ \text{elements}\}\crr
\ar=\ar \int_{[0,1]} x^{k-2}(1 - x)^{n-k}\Lambda(\dd x).
\eeqnn
Here $\Lambda$ is a finite measure on $[0, 1]$. Then the number of blocks jumps from $n$ to $n-k+1$ at rate ${n \choose k}\lambda_{n,k}.$ We refer the reader to \cite{P99} for details and additional analysis. We always assume $\Lambda(\{1\}) = 0$ so that it is impossible for all the blocks to coagulate simultaneously. Without lose of generality, to simplify the computation we also assume that $\Lambda$ is absolutely continuous with respect to the Lebesgue measure.

Let the splitting measure $\mu$ be the image of $\mu_\text{Frag}$ by the map $\pi \mapsto \sharp \pi - 1$. Since $\mu(\infty) = 0$, thereby ensuring that no block can be fragmented into infinitely many sub-blocks. The fragmentation is an opposite mechanism to coalescent, which was introduced by Bertoin \cite{B01} first, see also Bertoin \cite{B02, B03}. Upon the arrival of an atom $(t, \pi^f, j)$ of $\text{PPP}_\text{F}$ with $k:= \sharp \pi^f - 1$, given $\sharp \Pi(t-) = n$, if $j \le n$, then the $j^{th}$-block is fragmentated into $k+1$ blocks.  Therefore, the number of blocks jumps from $n$ to $n + k$ at rate $n\mu(k)$.

The simple EFC processes combine the above two mechanism. Recall that $N =(N_t, t \ge 0)$ is the corresponding block counting process with unspecified initial value. Then $N$ is a Markov process on state space $\bar{\mbb{N}}_+$ satisfying the Feller property by Foucart and Zhou \cite[Theorem 2.3]{FZ23}. The paths of $N$ are not monotone anymore. We cannot immediately deduce from the dynamics above whether the boundary $\infty$ can be reached or not. Let $\tau_\infty^+ := \inf\{t > 0: N_{t-} = \infty\}.$ Then by \cite[Proposition 2.11]{F22}, the process $(N_t, t < \tau_\infty^+)$ started from $n$ is Markovian on $\mbb{N}_+$ with its generator $\mathcal{L}$ acting on
\beqnn
\mathcal{D} := \left\{g: \mbb{N}_+ \mapsto \mbb{R}; \forall\ n \in \mbb{N}_+, \sum_{k \in \mbb{N}_+}|g(n+k)|\mu(k) < \infty\right\}
\eeqnn
as follows: for $n \in \mbb{N}_+$,
\beqlb\label{L}
\mathcal{L}g(n) :=\mathcal{L}^cg(n)+\mathcal{L}^fg(n)
\eeqlb
with
\beqnn
\mathcal{L}^cg(n):=\sum_{k=2}^n{n\choose k}\lambda_{n,k}[g(n-k+1)-g(n)]
\eeqnn
and
\beqnn
\mathcal{L}^fg(n):= n\sum_{k=1}^\infty\mu(k)[g(n+k)-g(n)].
\eeqnn
Here $\mathcal{L}^cg(n)$ vanishes if $n = 1.$ Then the density matrix $(q_{ij})_{i,j\in\mbb{N}_+}$ of $(N_t, t < \tau_\infty^+)$ is
\beqnn
q_{ij} =
\begin{cases}
i\mu(j-i),&  j \ge i+1, i \ge 1;\\
-\mu(\mbb{N}_+) - \sum_{k = 2}^i {i \choose k}\lambda_{i,k}, &   j = i \ge 1;\\
{i \choose k}\lambda_{i,k},&   j = i - k + 1, i \ge 1.
\end{cases}
\eeqnn 
Note that $\mbb{N}_+$ forms a communication class for the process $(N_t, t < \tau_\infty^+)$ when started from $n$ (see, e.g., \cite[Subsection 2.1]{FZ22}). Consequently, the process $(N_t, t < \tau_\infty^+)$ with initial value $n$ constitutes an  irreducible Markov process on $\mbb{N}_+$, governed by the generator $\mathcal{L}$ defined in \eqref{L}.

\subsection{Main Results}

In this paper, we are mainly interested in the critical regime with $d = \frac{b}{1 + \alpha}$ and $\beta = 1 + \alpha$ in Theorem \ref{FZ3.9}.
We first obtain the following estimate on $\Phi_\mu$ under Condition \ref{C1}.

\begin{proposition}\label{t10106}
Suppose that Condition \ref{C1} holds. Then $\Phi_\mu(n) \underset{n \rightarrow \infty}{\sim} \frac{b}{\alpha + 1}n(\log n)^{\alpha +1}$.
Therefore, the critical regime in Theorem \ref{FZ3.9} corresponds to $ \Phi_\Lambda(n) \underset{n \rightarrow \infty}{\sim} \Phi_\mu(n)$.
\end{proposition}

\proof
By Condition \ref{C1}, we have $\mu(k) = b(\log k)^\alpha k^{-2} + \epsilon(k)$ with $\epsilon(k)= o((\log k)^\alpha k^{-2})$ as $k \rightarrow \infty$. Then as $n\to\infty$,
\beqnn
\frac{\Phi_\mu(n)}{n}\ar=\ar \sum_{k = 1}^n\mu(k)k = \sum_{k = 1}^n\left(\frac{b(\log k)^\alpha}{k} + \epsilon(k)k\right)\cr
\ar=\ar b\int_1^n\frac{(\log x)^\alpha}{x}\dd x +  \int_1^n\epsilon(x)x \dd x  + O(1)\cr
\ar=\ar \frac{b}{\alpha + 1}(\log n)^{\alpha + 1} +  \int_1^n\epsilon(x)x \dd x  + O(1).
\eeqnn
Notice that $
\lim_{n \rightarrow \infty} \frac{\int_1^n\epsilon(x)x \dd x}{\int_1^n\frac{(\log x)^\alpha}{x}\dd x} = 0.$
Then
\beqnn
\frac{\Phi_\mu(n)}{n} = \frac{b}{\alpha + 1}(\log n)^{\alpha + 1} + o\left( (\log n)^{\alpha + 1}\right)
\eeqnn
as $n\to\infty$. The result of the proposition follows.
\qed

\begin{remark}
Observe from Proposition \ref{t10106} that $\sum_{k = 1}^{\lfloor rn\rfloor}\mu(k)k$ for any $r>0$ is of the same order of $\Phi_\mu(n)$. Very loosely put, $\Phi_\mu(n)$ represents the rate of new blocks produced in a fragmentation event when there are $n$ blocks in the simple EFC process.
In the critical case where
$$\Phi_\Lambda(n) \underset{n \rightarrow \infty}{\sim} \frac{b}{\alpha + 1}n(\log n)^{\alpha +1},$$
the boundary behavior then depends on the order of $\Phi_\Lambda(n)-\Phi_\mu(n)$.
Therefore, in the following Theorems, conditions are imposed on the difference between
$\Phi_\Lambda(n)$ and $\Phi_\mu(n)$  under Condition \ref{C1} on $\Phi_\mu$.
But to obtain more general results,  assumption is not always imposed on  $\Phi_\Lambda$.
\end{remark}

We first consider the case  $\alpha \in (0, 1]$.

\begin{theorem}\label{p3}
Assume that Condition \ref{C1} holds for $\alpha \in (0, 1]$. If there exists $m \in \mbb{N}_+$ such that
\beqlb\label{eq0120}
\liminf_{n \rightarrow \infty}\frac{ \Phi_\Lambda(n)-\Phi_\mu(n)}{\prod_{\ell = 0}^m \log^{(\ell)} n } > -\infty,
\eeqlb 
then the process $N$ does not explode.
\end{theorem}

\begin{theorem}\label{pp0406}
Suppose that Condition \ref{C1} holds for $\alpha \in (0,1]$. If
\beqlb\label{Phi_0518}
\liminf_{n \rightarrow \infty}\frac{\Phi_\Lambda(n) - \Phi_\mu(n)}{n\log n(\log\log n)^2} = \infty,
\eeqlb
then the process $N$ comes down from infinity.
\end{theorem}

In the above theorems, for $\alpha \in (0,1]$ we provide the conditions under which the process does not explode (Theorem \ref{p3}), as well as the those under which the process comes down from infinity (Theorem \ref{pp0406}).  We now proceed to provide the conditions for the process to stay infinite and those for
$N$ to explode. To this end, we need a condition on
$\Phi_\Lambda(n)$  for technical considerations.

\begin{condition}\label{CC2}
$0< \liminf_{n \rightarrow \infty}\frac{\Phi_\Lambda(n)}{n(\log n)^{\beta}} \le \limsup_{n \rightarrow \infty}\frac{\Phi_\Lambda(n)}{n(\log n)^{\beta}} < \infty$ for $\beta > 1$.
\end{condition}
Recall that $\Lambda$ is absolutely continuous with respect to the Lebesgue measure. 
Then Condition \ref{C2} holds if and only if there exist constants  $C_1, C_2 > 0$ such that
\beqlb\label{Lambda}
C_1\left(\log \frac{1}{x}\right)^{\beta - 1}\dd x \le \Lambda(\dd x) \le C_2\left(\log \frac{1}{x}\right)^{\beta - 1}\dd x
\eeqlb
holds for small $x$.
The equivalence of Condition \ref{CC2} and \eqref{Lambda} follows from \cite[Section 2.2]{F22}.  Without loss of generality, to simplify the computation we replace Condition \ref{CC2} with the following:
\begin{condition}\label{C2}
There exist constants  $C_1, C_2 > 0$ such that \eqref{Lambda} holds for all $x \in (0, 1)$.
\end{condition} 
Now we are ready to present the results regarding the  staying-infinite and explosion behaviors for $N$, respectively.

\begin{theorem}\label{p4}
Assume that Conditions \ref{C1} and \ref{C2} hold for $\alpha \in (0, 1]$. For $C_2 > 0$ in \eqref{Lambda}, if
\beqlb\label{Phi_0126}
\limsup_{n \rightarrow \infty}\left[\frac{\Phi_\Lambda(n) -  \Phi_\mu(n)}{n\log n\log\log n } + \frac{2C_2}{\beta}\log\log n\right] < \infty,
\eeqlb
then the process $N$ stays infinite.
\end{theorem}

\begin{corollary}\label{c2.7}
Assume that Condition \ref{C1} holds for $\alpha \in (0, 1]$ and $\Lambda(\dd x) = b(\log \frac{1}{x})^{\alpha}\dd x$. Then $\Phi_\Lambda(n) \underset{n \rightarrow \infty}{\sim} \frac{b}{\alpha + 1}n (\log n)^{\alpha + 1}$. Moreover, if
\beqnn 
\limsup_{n \rightarrow \infty} \frac{\Phi_\Lambda(n) -  \Phi_\mu(n)}{n\log n(\log\log n)^2 }<-\frac{2b}{\alpha + 1},  
\eeqnn 
then the process $N$ stays infinite.   
\end{corollary}
\proof
Notice that $\int_0^1 \Lambda(\dd x) = \int_0^1 b(\log \frac{1}{x})^{\alpha}\dd x < \infty$, and
\beqnn
\int_x^1\dd y\int_y^1 u^{-2}\Lambda(\dd u) \underset{x \rightarrow 0+}{\sim} \frac{b}{\alpha + 1} \left(\log \frac{1}{x}\right)^{\alpha + 1}.
\eeqnn
Then by \cite[Section 2.2]{F22}, we have $\Phi_\Lambda(n) \underset{n \rightarrow \infty}{\sim} \frac{b}{\alpha + 1}n (\log n)^{\alpha + 1}$. The second result follows from Theorem \ref{p4} directly.
\qed

\begin{theorem}\label{p0206}
Assume that Conditions \ref{C1} and \ref{C2} hold for $\alpha \in (0, 1]$. If 	
\beqlb\label{Phi_0518a}
\limsup_{n \rightarrow \infty}\frac{ \Phi_\Lambda(n)-\Phi_\mu(n)}{n\log n(\log\log n)^2} = -\infty, 
\eeqlb
then the process $N$ explodes.
\end{theorem}

Combining  Theorems \ref{p3} and \ref{pp0406}, we find conditions for $\infty$ to be an entrance boundary. Similarly, from Corollary \ref{c2.7} and Theorem \ref{p0206} we find conditions for $\infty$ to be an exit boundary.

\begin{corollary}\label{tt3}
Suppose that Condition \ref{C1} holds for $\alpha \in (0,1]$.
\begin{itemize}
\item[(i)] If  \eqref{Phi_0518} holds,
then $\infty$ is an entrance boundary;
\item[(ii)] If Condition \ref{C2} and  \eqref{Phi_0518a} hold,
then $\infty$ is an exit boundary.
\end{itemize}
\end{corollary}

We next consider the case of $\alpha > 1$ and  identify conditions of explosion/non-explosion, coming-down-from-infinity/staying-infinite, respectively, for process $N$.


\begin{theorem}\label{c0527a}
Suppose that Condition \ref{C1} holds for $\alpha > 1$. If
\beqlb\label{Phi_0527}
\liminf_{n \rightarrow \infty}\frac{\Phi_\Lambda(n) - \Phi_\mu(n)}{n(\log n)^{\alpha}} > \left( 2\log 2 -\frac{1}{2}\right)b,
\eeqlb
then the process $N$ does not explode.
\end{theorem}

\begin{theorem}\label{c0527}
Suppose that Condition \ref{C1} holds for  $\alpha > 1$. If \eqref{Phi_0527} holds, then the process $N$ comes down from infinity.
\end{theorem}

\begin{theorem}\label{p5}
Suppose that Conditions \ref{C1} and \ref{C2} hold for $\alpha > 1$.  If
\beqlb\label{Phi_0518b}
\limsup_{n \rightarrow \infty}\frac{\Phi_\Lambda(n) - \Phi_\mu(n)}{n(\log n)^{\alpha} } < \left( 2\log 2 -\frac{1}{2}\right)b,
\eeqlb
then the process $N$ stays infinite.
\end{theorem}

\begin{theorem}\label{p0206a}
Suppose that Conditions \ref{C1} and \ref{C2} hold for $\alpha > 1$.	If \eqref{Phi_0518b} holds, then the process $N$ explodes.
\end{theorem}

The above results further allow to specify the type of boundary  at $\infty$.

\begin{corollary}\label{t0304}
Suppose that Condition \ref{C1} holds for $\alpha > 1$.
\begin{itemize}
\item[(i)]  If \eqref{Phi_0527} holds, then $\infty$ is an entrance boundary for $N$;
\item[(ii)] If Condition \ref{C2} and \eqref{Phi_0518b} hold,
then $\infty$ is an exit boundary for $N$.
\end{itemize}
\end{corollary}

Moreover, for the case of $\alpha > 1,$ we have the following result.

\begin{corollary}
If $\mu(n) = b(\log n)^\alpha n^{-2}$ and $\Phi_\Lambda(n) = \frac{b}{1+\alpha}n(\log n)^{1+\alpha}$ for $b > 0$ and $\alpha > 1$, then $\infty$ is an exit boundary.
\end{corollary}

\section{Preliminaries}

\subsection{Some known results}\label{s3}
We recall some previous results in this subsection, especially the coupling for EFC process established in \cite{F22}. Recall that $\Pi(0) = \{\Pi_1(0), \Pi_2(0), \cdots\}$ is proper. Without loss of generality, we assume $\sharp\Pi(0) = \infty$ almost surely.  For any $n \in \bar{\mbb{N}}_+$, a process $\Pi^{(n)}:=(\Pi^{(n)}(t), t \ge 0)$, started from $\Pi^{(n)}(0) := \{\Pi_1(0), \cdots, \Pi_n(0)\}$, can be constructed from $\text{PPP}_\text{C}$ and $\text{PPP}_\text{F}$, see Subsection \ref{ss21}. Then the process $\Pi^{(n)}$ follows all coagulations and fragmentations involving integers belong to $\cup_{i=1}^n \Pi_i(0).$ We refer to \cite[Section 3]{F22} for details of $\Pi^{(n)}$.  One sees that $\Pi^{(\infty)}= \Pi$ almost surely. Let $N^{(n)}:=(N^{(n)}_t, t \ge 0)$ be the block counting process of $\Pi^{(n)}$, i.e., $N^{(n)}_t : = \sharp \Pi^{(n)}(t)$ with $N^{(n)}_0 = n$. Then $N^{(\infty)}_t =\sharp \Pi^{(\infty)}(t) = \sharp \Pi(t)$ for all $t \ge 0$ almost surely with $N^{(\infty)}_0 = \infty$. The following result plays a crucial role in our proofs for explosion/non-explosion.

\begin{lemma}[Lemma 3.4 of \cite{F22}]\label{Fl3.4}
Assume that $\Pi(0)$ consists of blocks with infinite sizes. For any $n \in \mbb{N}_+,$ set $\tau_\infty^{n,+} := \inf\{t > 0: N_{t-}^{(n)} = \infty\}$. The process $(N_t^{(n)}, t < \tau_\infty^{n,+})$ has the same law as $(N_t, t < \tau_\infty^+)$ started from $n.$ Moreover almost surely, for all $n \in \mbb{N}_+$ and all $t \ge 0$,
\beqnn
N^{(n)}_t \le N^{(n+1)}_t\ \text{for\ all}\ n \ge 1\quad \text{and}\quad \lim_{n \rightarrow \infty} N^{(n)}_t = N_t^{(\infty)}.
\eeqnn
\end{lemma}

According to the above result and \cite[Proposition 2.11]{F22}, the process $(N^{(n)}_t, t < \tau_\infty^{n,+})$ is a Markov process on $\mbb{N}_+$ with initial value $n$ and generator $\mathcal{L}$ given by \eqref{L}.

Now we introduce a partition-valued process $(\Pi^m(t), t \ge 0)$, in which every fragmentation creates at most $m$ new blocks. For any $m \in \mbb{N}_+$, define a map
\beqnn
r_m: \pi \mapsto (\pi_1, \cdots, \pi_m, \cup_{k = m + 1}^\infty \pi_k),
\eeqnn
which maps $\mathcal{P}_\infty$ to partitions with at most $m+1$ blocks. Set $\mu_\text{Frag}^m := \mu_\text{Frag}\circ r_m^{-1}$. Then we construct $\Pi^m:=(\Pi^m(t), t \ge 0)$ with initial value $\Pi(0)$ through Poissonian construction similar to that for $(\Pi(t), t \ge 0)$ where the coagulation measure $\mu_{\text Coag}$ is kept the same but the fragmentation measure $\mu_\text{Frag}$ is replaced by $\mu_\text{Frag}^m$.
Similar to $(\Pi^{(n)})$, a monotone coupling $\Pi^{m,(n)}:=(\Pi^{m,(n)}(t), t \ge 0)$ to $\Pi^m$ can also be built via the same Poissonian construction with $\Pi$ replaced by $\Pi^m$. 
We refer to \cite[Subsection 3.2]{F22} for details on the construction of $\Pi^m$ and $\Pi^{m,(n)}$.

Set $N_m^{(n)}(t):= \sharp\Pi^{m,(n)}(t)$ for each $n \in \bar{\mbb{N}}_+$ and $t \ge 0.$ Then $N_m^{(\infty)}(t) = \sharp\Pi^{m,(\infty)}(t) = \sharp\Pi^{m}(t)$ for $t \ge 0$ almost surely. For any $m \in \mbb{N}_+$, let $\mu_m$ be the image of $\mu_\text{Frag}^m$ by the map $\pi \mapsto \sharp\pi - 1$, i.e., $\mu_m(l) = \mu(l)$ if $l \le m - 1$ and $\mu_m(m) = \sum_{k = m}^\infty\mu(k).$ Set the operator $\mathcal{L}^m$ acting on any function $g: \mbb{N}_+ \mapsto \mbb{R}_+$ as follows:
\beqlb\label{Lm}
\mathcal{L}^mg(n) := \mathcal{L}^cg(n) + n\sum_{l=1}^m\mu_m(l)[g(l+n)-g(n)] \quad \text{for\ all}\ n \in \mbb{N}_+.
\eeqlb
The following result follows from elementary calculation. We state it without proof.
\begin{lemma}\label{ll0406}
Let $g$ be a positive decreasing function on $\mbb{N}_+.$ Then for every $m \in \mathbb{N}_+$, $\mathcal{L}^m g(n) \ge \mathcal{L}g(n)$ holds for every $n \in \mathbb{N}_+$.
\end{lemma}

The following lemma
guarantees that the process $(N_t^{(n)}, t \ge 0)$ can be approximated by a non-decreasing sequence of non-explosive processes. It plays a key role in establishing the entrance boundary.

\begin{lemma}[Lemma 3.8 of \cite{F22}]\label{Fl3.8}
Process $(N_m^{(n)}(t), t \ge 0)$ is a nonexplosive Markov process started from $n$ with generator $\mathcal{L}^m$ given by \eqref{Lm}. Moreover, almost surely for any $m, n \in \mbb{N}_+$ and all $t \ge 0$, $N_m^{(n)}(t) \le N_m^{(n+1)}(t),\  N_m^{(\infty)}(t) \le N_{m+1}^{(\infty)}(t)$
and
\beqnn
\lim_{m \rightarrow \infty}N_m^{(\infty)}(t) = N_t^{(\infty)}.
\eeqnn
\end{lemma}

For any $m,a \in \mathbb{N}_+$ and $n \in \bar{\mbb{N}}_+$ with $n > a$, we consider the first entrance times $\tau_{a,m}^{n,-} := \inf\{t > 0: N_m^{(n)}(t) \le a\}$ and $\tau_{a}^{n,-} := \inf\{t > 0: N_t^{(n)} \le a\}$. We simplify them as $\tau_{a,m}^-$, $\tau_a^-$, respectively, if the initial value $n$ is obvious.
\begin{lemma}[Lemma 3.10 of \cite{F22}]\label{Fl3.10}
For any $a, m\in \mbb{N}_+$, $\lim_{n \rightarrow \infty}\tau_{a,m}^{n,-} = \tau_{a,m}^{\infty,-}$ almost surely. If moreover for any $m \in \mbb{N}_+$ and any $s > 0$, $N_m^{(\infty)}(s) < \infty$, then $\lim_{m \rightarrow \infty}\tau_{a,m}^{\infty,-} = \tau_a^{\infty,-}$ almost surely.
\end{lemma}

\subsection{Estimations of $\Lambda$ and $\mu$}

In this subsection we present estimates concerning the coagulation and fragmentation rates.
\begin{proposition}\label{Phi2}
For any $n\geq 2$, we have
\beqnn
\sum_{k = 2}^n \binom{n}{k}\lambda_{n, k}(k - 1)^2 = \Lambda[0,1)n(n-1) - \Phi_\Lambda(n).
\eeqnn
Furthermore, $n^{-2}\sum_{k = 2}^n \binom{n}{k}\lambda_{n, k}(k - 1)^2 = O(1)$ as $n\to\infty$.
\end{proposition}
\proof
By the definition of $\lambda_{n,k}$ and $(k-1)^2 = k(k-1) - (k - 1),$ it is easy to see that
\beqnn
0 < \sum_{k = 2}^n \binom{n}{k}\lambda_{n, k}(k - 1)^2 \ar=\ar \sum_{k = 2}^n  \binom{n}{k}\lambda_{n, k} k(k - 1) - \Phi_\Lambda(n)\crr
\ar=\ar \int_{[0,1)} \sum_{k = 2}^n \frac{n!}{(k - 2)!(n - k)!}x^{k - 2}(1 - x)^{n - k}\Lambda(\dd x) - \Phi_\Lambda(n)\crr
\ar=\ar \Lambda[0,1)n(n-1) - \Phi_\Lambda(n).
\eeqnn
The result follows.
\qed

\begin{proposition}\label{p0217}
As $n\to\infty$,
\beqnn
\sum_{k=2}^n{n\choose k}\lambda_{n,k} \log \left(1 - \frac{k - 1}{n}\right) = -\frac{\Phi_\Lambda(n)}{n} + O(1).
\eeqnn
\end{proposition}
\proof
For a fixed constant $C > 0$, we have
\beqnn
\log \left(1 - \frac{C}{n}\right) = - \frac{C}{n} - \frac{C^2}{2n^2} - C^2 \epsilon_n,
\eeqnn
where $\epsilon_n = o(n^{-2})$ as $n \rightarrow \infty$.
By the above and Proposition \ref{Phi2}, we have
\beqnn
\ar\ar\sum_{k=2}^n{n\choose k}\lambda_{n,k} \log \left(1 - \frac{k - 1}{n}\right)\crr \ar\ar\quad= -\sum_{k=2}^n{n\choose k}\lambda_{n,k} \left(\frac{k - 1}{n} + \frac{(k - 1)^2}{2n^2} +  \epsilon_n(k - 1)^2 \right)\crr
\ar\ar\quad= -\frac{\Phi_\Lambda(n)}{n} - \frac{1}{2n^2}\sum_{k=2}^n{n\choose k}\lambda_{n,k}(k - 1)^2  - \epsilon_n \sum_{k=2}^n{n\choose k}\lambda_{n,k}(k - 1)^2\cr
\ar\ar\quad= -\frac{\Phi_\Lambda(n)}{n} + O(1).
\eeqnn
The result follows.
\qed

\begin{proposition}\label{l0225}
Assume that Condition \ref{C2} holds. Then for any $\delta \in (0, 1)$, as $n \rightarrow \infty,$ we have
\beqnn
\sum_{k = \lfloor\delta n\rfloor + 1}^n {n \choose k}\lambda_{n,k} \le \frac{C_2}{\beta}\delta^{-\beta}  + O(n^{-1}),
\eeqnn
where $C_2$ is the constant given in \eqref{Lambda}.
\end{proposition}
\proof
Recall that $\beta > 1$ and Condition \ref{C2} holds if and only if \eqref{Lambda} holds. Take $n > \beta/\delta$. For any $k \in [\lfloor\delta n\rfloor + 1, n]$, by \eqref{Lambda} and  the fact of $\log (x^{-1}) \le x^{-1}$ for any $x \in (0, 1)$, it follows that
\beqnn
\lambda_{n,k}\ar=\ar \int_{[0,1)} x^{k-2}(1-x)^{n-k}\Lambda(\dd x)\cr
\ar\le\ar C_2\int_0^1 x^{k-2}(1-x)^{n-k}\left(\log \frac{1}{x}\right)^{\beta - 1}\dd x\cr
\ar\le\ar C_2\int_0^1 x^{k-2}(1-x)^{n-k}x^{-(\beta - 1)} \dd x \cr
\ar=\ar C_2\text{ Beta}(k-\beta, n-k+1).
\eeqnn
We further obtain that
\beqnn
{n\choose k}\lambda_{n,k} \ar\le\ar C_2{n\choose k}\text{ Beta}(k-\beta, n-k+1) \cr
\ar=\ar C_2{n \choose k}\frac{\Gamma(k-\beta)\Gamma(n-k+1)}{\Gamma(n -\beta+1)} = C_2\frac{n!}{\Gamma(n -\beta+1)}\frac{\Gamma(k-\beta)}{k!}\cr
\ar=\ar C_2\frac{n\cdots (n - \lfloor \beta \rfloor + 1)\Gamma(n - \lfloor \beta \rfloor + 1)}{\Gamma(n -\beta+1)}\frac{\Gamma(k-\beta)}{k(k-1)\cdots (k - \lfloor \beta\rfloor)\Gamma(k - \lfloor \beta\rfloor)}.
\eeqnn
By Gautschi's inequality, i.e., for any $x>0$ and $s\in (0,1)$,
\beqlb\label{Gautschi}
x^{1-s}<\frac{\Gamma(x+1)}{\Gamma(x+s)}<(x+1)^{1-s},
\eeqlb
one sees that
\beqnn
(n - \lfloor\beta\rfloor)^{\beta  - \lfloor\beta\rfloor}\le \frac{\Gamma(n - \lfloor\beta\rfloor + 1)}{\Gamma(n -\beta+1)} \le (n - \lfloor\beta\rfloor + 1)^{\beta  - \lfloor\beta\rfloor}
\eeqnn
and
\beqnn
(k - \lfloor\beta\rfloor - 1)^{\beta - \lfloor\beta\rfloor} \le \frac{\Gamma(k - \lfloor\beta\rfloor)}{\Gamma(k-\beta)} \le (k - \lfloor\beta\rfloor)^{\beta - \lfloor\beta\rfloor}.
\eeqnn
It implies that
\beqnn
\frac{n!}{\Gamma(n -\beta+1)} \le n^\beta \quad \text{and}\quad  \frac{\Gamma(k-\beta)}{k!} \le (k-\lfloor\beta\rfloor-1)^{-\beta - 1}.
\eeqnn
Then we have
\beqnn
\sum_{k = \lfloor\delta n\rfloor+1}^n {n\choose k}\lambda_{n,k} \le C_2 n^\beta\sum_{k = \lfloor\delta n\rfloor+1}^n (k-\lfloor\beta\rfloor-1)^{-\beta - 1}.
\eeqnn
The result follows.
\qed

\begin{proposition}\label{p0211}
Suppose that Condition \ref{C1} holds. Then as $n\to\infty$,
\beqnn
n\sum_{k = 1}^n \mu(k)\log\left(1 + \frac{k}{n}\right) = \frac{\Phi_\mu(n)}{n} - \frac{b(\log n)^\alpha}{2} + o((\log n)^\alpha).
\eeqnn

\end{proposition}
\proof
For a fixed constant $C > 0$, notice that
\beqnn
\log\left(1 + \frac{C}{n}\right) = \frac{C}{n} - \frac{C^2}{2n^2} + C^2 \epsilon_n
\eeqnn
with $\epsilon_n = o(n^{-2})$ as $n \rightarrow \infty$, then we have
\beqlb\label{eq0218a}
n\sum_{k = 1}^n \mu(k)\log\left(1 + \frac{k}{n}\right) \ar=\ar 	n\sum_{k = 1}^n \mu(k)\left(\frac{k}{n} - \frac{k^2}{2n^2} + \epsilon_n k^2\right)\cr
\ar=\ar  \frac{\Phi_\mu(n)}{n} - \frac{1}{2n}\sum_{k = 1}^n \mu(k)k^2  + n\epsilon_n \sum_{k = 1}^n \mu(k)k^2.
\eeqlb
By Condition \ref{C1}, using the notation in the proof of Proposition \ref{t10106}, as $n\to\infty$, it follows that
\beqlb\label{eq310}
\sum_{k = 1}^n \mu(k)k^2 \ar=\ar \sum_{k = 1}^n \left(b(\log k)^\alpha + \epsilon(k)k^2\right)\cr
\ar=\ar b\int_1^n (\log x)^\alpha \dd x + \int_1^n \epsilon(x)x^2 \dd x + O(1)\cr
\ar=\ar b\int_1^n (\log x)^\alpha \dd x + o\left(\int_1^n (\log x)^\alpha \dd x \right) + O(1).
\eeqlb
Notice that
\beqlb\label{eq11}
(1 - a x)^\alpha = 1 - a \tilde{\epsilon}(x)
\eeqlb
for $x \in (0, 1/a)$ with $a > 0$ being a fixed constant, where $\tilde{\epsilon}(x) = o(x)$ as $x \rightarrow 0+$. By the change of variable $y = x/n$, as $n \rightarrow \infty,$ we have
\beqlb\label{eq310a}
\int_1^n (\log x)^\alpha \dd x \ar=\ar n(\log n)^\alpha \int_{1/n}^1 \left(1 + \frac{\log y}{\log n}\right)^\alpha \dd y\cr
\ar=\ar n(\log n)^\alpha \left(1 - \frac{1}{n}\right) + n(\log n)^\alpha \tilde{\epsilon} \left(\frac{1}{\log n}\right) \int_{1/n}^1 \log y  \dd y\cr
\ar=\ar  n(\log n)^\alpha \left(1 - \frac{1}{n}\right) +  n(\log n)^{\alpha } \tilde{\epsilon} \left(\frac{1}{\log n}\right)\left[\frac{\log n}{n} - \left(1 - \frac{1}{n}\right) \right]\cr
\ar=\ar n(\log n)^\alpha + o(n(\log n)^\alpha).
\eeqlb
The second equality follows from \eqref{eq11} by taking $a = -\log y$ and $x = \frac{1}{\log n}$. By \eqref{eq310} and \eqref{eq310a}, as $n \rightarrow \infty$, one obtains
\beqlb\label{mu2}
\sum_{k = 1}^n \mu(k)k^2 = bn(\log n)^\alpha + o(n(\log n)^\alpha).
\eeqlb
Then the result follows from \eqref{eq0218a} and \eqref{mu2}.
\qed

\begin{proposition}\label{mu_tail}
Suppose that Condition \ref{C1} holds. Then as $n\to\infty$,
\beqnn
n\sum_{k = n + 1}^\infty \mu(k)\log\left(1 + \frac{k}{n}\right)=  (2\log 2)b(\log n)^\alpha + o((\log n)^{\alpha}).
\eeqnn
\end{proposition}
\proof By integration by parts and change of variable $y = x/n$,
as $n\to\infty$, we have
\beqlb\label{1212}
\ar\ar n\sum_{k=n+1}^\infty \frac{(\log k)^\alpha}{k^2}\log\left(1+\frac{k}{n}\right) \cr
\ar\ar\quad= n\int_n^\infty \frac{(\log x)^\alpha}{x^2}\log\left(1+\frac{x}{n}\right)\dd x + O(1)\cr
\ar\ar\quad=-n\int_n^\infty (\log x)^\alpha\log\left(1+\frac{x}{n}\right)\dd\frac{1}{x} + O(1)\cr
\ar\ar\quad=(\log n)^\alpha\log 2+n\int_n^\infty \frac{1}{x}\left(\frac{\alpha(\log x)^{\alpha-1}}{x}\log\left(1+\frac{x}{n}\right)+ \frac{(\log x)^\alpha}{1+x/n}\frac{1}{n}\right)\dd x + O(1)\cr
\ar\ar\quad=(\log n)^\alpha\log 2+\int_1^\infty \frac{\alpha(\log y+\log n)^{\alpha-1}}{y^2}\log(1+y)\dd y\cr
\ar\ar\qquad+\int_1^\infty \frac{(\log y+\log n)^\alpha}{(1+y)y}\dd y + O(1)\cr
\ar\ar\quad=: (\log n)^\alpha\log 2+ I_n^1 + I_n^2 + O(1).
\eeqlb
Notice that $I_n^1$ is finite when $\alpha \in (0, 1]$. For the case of $\alpha > 1$, as $n\to\infty$, we have
\beqnn
0\le I_n^1
\ar=\ar \int_1^n \frac{\alpha(\log y+\log n)^{\alpha-1}}{y^2}\log(1+y)dy\cr
\ar\ar  + \int_{n}^\infty \frac{\alpha(\log y+\log n)^{\alpha-1}}{y^2}\log(1+y)\dd y\cr
\ar\le\ar \alpha(2\log n)^{\alpha - 1} \int_1^n\frac{\log(1 + y)}{y^2}\dd y + \int_n^\infty\frac{\alpha(2\log y)^{\alpha-1}\log(1 + y)}{y^2}\dd y\cr
\ar=\ar O((\log n)^{\alpha - 1}).
\eeqnn
Then for $\alpha > 0,$
\beqlb\label{eqI1}
I_n^1 \ar=\ar o((\log n)^{\alpha}) \qquad \text{as}\ n \rightarrow \infty.
\eeqlb
On the other hand, there exists a constant $C_\alpha > 0$ depending on $\alpha$ such that
\beqnn
(1 + x)^\alpha \le 1 + C_\alpha x, \,\,\, x \in (0, 1).
\eeqnn
Then for the third term of \eqref{1212}, as $n \rightarrow \infty$, we have
\beqnn
I_n^2
\ar=\ar \int_1^n \frac{(\log y+\log n)^\alpha}{(1+y)y}\dd y+\int_n^\infty \frac{(\log y+\log n)^\alpha}{(1+y)y}\dd y\cr
\ar\le\ar (\log n)^\alpha\int_1^n \frac{(1 + \frac{\log y}{\log n})^\alpha}{(1+y)y}dy+\int_n^\infty \frac{(2\log y)^\alpha}{(1+y)y}\dd y\cr
\ar\le\ar  (\log n)^\alpha\int_1^n \frac{1}{(1+y)y}dy +C_\alpha(\log n)^{\alpha - 1}\int_1^n \frac{\log y}{(1+y)y}dy + O(1)\crr
\ar\le\ar (\log 2)(\log n)^\alpha +C_\alpha(\log n)^{\alpha - 1}\int_1^n \frac{\log y}{y^2}dy + O(1)\crr
\ar=\ar (\log 2)(\log n)^\alpha + o((\log n)^\alpha)
\eeqnn
and
\beqnn
I_n^2 \ar\ge\ar (\log n)^\alpha\int_1^n\frac{1}{y(1+y)}\dd y\cr
\ar\ge\ar (\log 2)(\log n)^\alpha - \frac{(\log n)^\alpha}{n},
\eeqnn
which implies
\beqlb\label{eqI2}
I_n^2 = (\log 2)(\log n)^\alpha + o((\log n)^\alpha), \qquad \text{as}\ n \rightarrow \infty.
\eeqlb
Then it follows from \eqref{1212}, \eqref{eqI1} and \eqref{eqI2} that
\beqlb\label{eq0218b}
n\sum_{k=n+1}^\infty \frac{(\log k)^\alpha}{k^2}\log\left(1+\frac{k}{n}\right)  = (2\log 2)(\log n)^\alpha + o((\log n)^\alpha)
\eeqlb
as $n \rightarrow \infty.$
Moreover, recall that $\mu(k) = b(\log k)^\alpha k^{-2} + \epsilon(k)$ with $\epsilon(k)=o((\log k)^\alpha /k^{2})$ as $k \rightarrow \infty$. Then we have
\beqnn
\lim_{n \rightarrow \infty}\frac{n\sum_{k=n+1}^\infty \epsilon(k)\log\left(1+\frac{k}{n}\right)}{n\sum_{k=n+1}^\infty \frac{(\log k)^\alpha}{k^2}\log\left(1+\frac{k}{n}\right)} \ar=\ar \lim_{n \rightarrow \infty} \frac{\int_1^\infty \epsilon(nx)\log(1+x)\dd x}{\int_1^\infty \frac{(\log(nx))^\alpha}{(nx)^2} \log (1+x)\dd x}\crr
\ar=\ar \lim_{n \rightarrow \infty} \frac{\int_1^\infty \frac{\epsilon(nx)}{(\log n)^\alpha n^{-2}}\log(1+x)\dd x}{\int_1^\infty \frac{(1 + \frac{\log x}{\log n})^\alpha}{x^2} \log (1+x)\dd x} = 0,
\eeqnn
which, by \eqref{eq0218b}, implies that
\beqnn
n\sum_{k=n+1}^\infty \epsilon(k)\log\left(1+\frac{k}{n}\right) = o((\log n)^\alpha), \qquad \text{as}\ n \rightarrow \infty.
\eeqnn
It follows that
\beqnn
n\sum_{k = n + 1}^\infty \mu(k)\log\left(1 + \frac{k}{n}\right) \ar=\ar  n\sum_{k=n+1}^\infty \left(\frac{b(\log k)^\alpha}{k^2} + \epsilon(k)\right)\log\left(1+\frac{k}{n}\right)\cr
\ar=\ar (2\log 2)b(\log n)^\alpha + o((\log n)^{\alpha})
\eeqnn
as $n \rightarrow \infty.$
This completes the proof.
\qed

The next result is similar to Proposition \ref{mu_tail} whose proof is omitted.
\begin{proposition}\label{c_mu_tail}
Suppose that Condition \ref{C1} holds. Then as $n\to\infty$,
\beqnn
n\sum_{k = n+1}^\infty\mu(k)\left(\log\left(1+\frac{k}{n}\right)\right)^2 = O((\log n)^\alpha).
\eeqnn
\end{proposition}

\section{Some criteria for boundary classification}

\subsection{Criteria for Markov processes on $\mbb{N}_+$}\label{ss4.1}

In this subsection, we present some criteria for Markov processes on $\mbb{N}_+$, adapting techniques pioneered by Mufa Chen to classify boundaries for Markov jump processes through conditions on generators (see \cite{C86a, C86b} and Theorems 2.25–2.27 of \cite{C04}, along with some developments in \cite{MT93}). These methods have been extended to study boundary behaviors of SDEs linked to continuous-state branching processes in \cite{LYZ19, MYZ21, RXYZ22, MZ23}.

Let $Z := (Z_t, t \ge 0)$ be an irreducible continuous-time Markov chain on $\mbb{N}_+$ with generator $\mcr{L}.$ For any positive constants $a, b$ with $a < n < b$ and $Z_0 = n$, let $\kappa_a^{n,-} := \inf\{t > 0: Z_t \le a\}$, $\kappa_b^{n,+} := \inf\{t > 0: Z_t \ge b\}$, $\kappa_{a,b}^n := \kappa_a^{n,-}\wedge\kappa_b^{n,+}$ and $\kappa_\infty^{n,+} := \lim_{b \rightarrow \infty}\kappa_b^{n,+}.$ For simplicity, we adopt the notations $\kappa_a^-$, $\kappa_b^+$, $\kappa_{a,b}$ and $\kappa_\infty^+$ when the initial value is specified. The following is a classical sufficient condition for non-explosion. We refer to, for instance,  Chow and Khasminskii \cite{CK11} and the references therein.
\begin{proposition}\label{p0607}
Let $g$ be a non-decreasing function on $\mbb{N}_+$ satisfying $\lim_{n \rightarrow \infty}g(n) = \infty.$ If there exists a constant $C > 0$ such that
\beqlb\label{F-L}
\mcr{L}g(n) \le Cg(n),\quad n \ge 1,
\eeqlb
then the process $Z$ does not explode.
\end{proposition}

By an adaption of \cite[Lemma 3.3]{MYZ21}, we have the following result.
\begin{proposition}\label{CDI}
Suppose that there exist a bounded  strictly positive function $g(u)$ on $\mathbb{N}_+$ with $\limsup_{u\to\infty} g(u)>0$,
and an eventually strictly positive function $d(a)$ on $\mathbb{N}_+$ with $\lim_{a\to\infty}d(a)=\infty$
such that for any $a\in \mathbb{N}_+$,
\beqlb\label{L0208}
\mcr{L}g(u) \ge d(a)g(u)
\eeqlb
holds for all $u>a $.
\begin{itemize}
\item If process $Z$ does not explode, then for any $t > 0$,
\beqlb\label{eq0306a}
\lim_{a\to\infty}\lim_{n\to\infty}\mbf{P}_n\{\kappa_a^-<t\}\ge \frac{\limsup_{u \rightarrow \infty}g(u)}{\sup_u g(u)};
\eeqlb
\item If $\lim_{n \rightarrow \infty}\mbf{P}_n\{{ \kappa_a^-} < t \wedge \kappa^+_\infty\} = 0$ for any $a, t > 0$, then it explodes.
\end{itemize}
\end{proposition}

\proof
For all large $0<a<n <b$, we have
\beqnn
\mbf{E}_n\big[g(Z_{t_0\wedge\kappa_{a,b}})\big]
\ar=\ar
g(n)
+\mbf{E}_n\Big[\int_0^{t_0\wedge\kappa_{a,b}}\mcr{L}g(Z_s)\dd s\Big]\cr
\ar=\ar
g(n)
+\int_0^{t_0} \mbf{E}_n\Big[\mcr{L}g(Z_s)1_{\{s\le \kappa_{a,b}\}}\Big] \dd s
\eeqnn
and then by integration by parts,
\beqnn
\ar\ar
\mbf{E}_n\big[g(Z_{t_0\wedge\kappa_{a,b}})\big]\e^{-d(a)t_0} \cr
\ar\ar\quad=
g(n)+\int_0^{t_0} \mbf{E}_n\big[g(Z_{s\wedge\kappa_{a,b}})\big]\dd (\e^{-d(a)s})
+\int_0^{t_0} \e^{-d(a)s}\dd \big(\mbf{E}_n\big[g(Z_{s\wedge\kappa_{a,b}})\big]\big) \cr
\ar\ar\quad=
g(n)
-d(a)\int_0^{t_0}\mbf{E}_n\big[g(Z_{s\wedge\kappa_{a,b}})\e^{-d(a)s}\big]\dd s
+\int_0^{t_0} \e^{-d(a)s}\mbf{E}_n\big[\mcr{L}g(Z_s)1_{\{s\le\kappa_{a,b}\}}\big]\dd s \cr
\ar\ar\quad\ge
g(n)
-d(a)\int_0^{t_0}\mbf{E}_n\big[g(Z_{s\wedge\kappa_{a,b}})\big]\e^{-d(a)s}\dd s
+d(a)\int_0^{t_0} \e^{-d(a)s}\mbf{E}_n\big[g(Z_s)1_{\{s\le\kappa_{a,b}\}}\big]\dd s,
\eeqnn
where the above inequality is from \eqref{L0208}. It implies that
\beqnn
g(n)\le\mbf{E}_n\big[g(Z_{t_0\wedge\kappa_{a,b}})\e^{-d(a)t_0}\big]
+d(a)\mbf{E}_n\Big[\int_0^{t_0}g(Z_{\kappa_{a,b}})\e^{-d(a)s}
1_{\{s>\kappa_{a,b}\}}\dd s\Big].
\eeqnn
Letting $t_0\to\infty$ in the above inequality and using the dominated convergence we obtain
\beqnn
g(n)
\ar\le\ar
d(a)\mbf{E}_n\Big[g(Z_{\kappa_{a,b}})\int_{\kappa_{a,b}}^\infty\e^{-d(a)s}\dd s\Big]
=\mbf{E}_n\big[g(Z_{\kappa_{a,b}}) \e^{-d(a)\kappa_{a,b}}\big]\cr
\ar=\ar
\mbf{E}_n\Big[g(Z_{\kappa_{a,b}})
\e^{-(\kappa_a^-\wedge\kappa^+_b)d(a)}
\big(1_{\{\kappa_b^+<\kappa^-_a\}}+1_{\{\kappa_a^-<t \wedge\kappa^+_b\}}
+1_{\{t \le\kappa^-_a<\kappa_b^+\}}\big)\Big]
\eeqnn
for any $t > 0$. It follows that
\beqnn
g(n)\ar\leq\ar
\limsup_{b\to\infty} \mbf{E}_n\Big[  g(Z_{\kappa_{a,b}})
\e^{-(\kappa_a^-\wedge\kappa^+_b)d(a)}
\big(1_{\{\kappa_b^+<\kappa^-_a\}}+1_{\{\kappa_a^-<t \wedge\kappa^+_b\}}
+1_{\{t \le\kappa^-_a<\kappa_b^+\}}\big)\Big] \cr
\ar\leq\ar
\mbf{E}_n\Big[\limsup_{b\to\infty}  g(Z_{\kappa_{a,b}})
\e^{-(\kappa_a^-\wedge\kappa^+_\infty)d(a)}
\big(1_{\{\kappa_\infty^+<\tau^-_a\}}+1_{\{\kappa_a^-<t \wedge\kappa^+_\infty\}}
+1_{\{t \le\kappa^-_a<\kappa_\infty^+\}}\big)\Big] \cr
\ar\le\ar
\limsup_{u\to\infty}g(u)
\mbf{P}_n\{\kappa_\infty^+<\infty\}
+\mbf{E}_n[g(Z_{\kappa_a^-}) 1_{\{\kappa_a^-<t \wedge\kappa^+_\infty \}}]\crr
\ar\ar +\e^{-d(a) t }\mbf{E}_n[g(Z_{\kappa_a^-})1_{\{\kappa^-_a<\kappa^+_\infty\}} ] \cr
\ar\le\ar \limsup_{u\to\infty}g(u)\mbf{P}_n \{\kappa_\infty^+<\infty\}
+\sup_u g(u)\left[\mbf{P}_n\{\kappa_a^-< t\wedge \kappa^+_\infty \}+\e^{-d(a) t }\right].
\eeqnn
Letting $n\to\infty$, we have
\beqlb\label{F_L_a}
\limsup_{n\to\infty}g(n)
\ar\leq\ar \limsup_{u\to\infty}g(u)
\limsup_{n\to\infty}\mbf{P}_n\{\kappa_\infty^+<\infty\}\cr
\ar\ar +\sup_{u}g(u)\limsup_{n\to\infty}\mbf{P}_n\{\kappa_a^-<t \wedge \kappa_\infty^+\}+\sup_{u}g(u)\e^{-d(a) t}\cr
\ar=\ar \limsup_{u\to\infty}g(u)
\lim_{n\to\infty}\mbf{P}_n\{\kappa_\infty^+<\infty\}\cr
\ar\ar +\sup_{u}g(u)\lim_{n\to\infty}\mbf{P}_n\{\kappa_a^-<t\wedge\kappa_\infty^+ \}+\sup_{u}g(u)\e^{-d(a) t}.
\eeqlb

If there is no explosion, i.e., $\lim_{n \rightarrow \infty}\mbf{P}_n\{\kappa_\infty^+ < \infty\} = 0,$ by letting $a\to\infty$ we have
\beqnn
\limsup_{n\to\infty}g(n) \ar\le\ar \sup_{u}g(u)\limsup_{a\to\infty}
\lim_{n\to\infty}\mbf{P}_n\{\kappa_a^-<t \wedge\kappa_\infty^+\}\cr
\ar=\ar \sup_{u}g(u)\limsup_{a\to\infty}
\lim_{n\to\infty}\mbf{P}_n\{\kappa_a^-<t \}.
\eeqnn
Observing that $\lim_{n\to\infty}\mbf{P}_n\{\kappa_a^-<t \} $ is increasing in $a$, then for any $t >0$, \eqref{eq0306a} holds.

Now we consider the case of $\lim_{n \rightarrow \infty}\mbf{P}_n\{\kappa_a^- < t \wedge \kappa_\infty^+\} = 0$ for any $a, t > 0$. By \eqref{F_L_a} and letting $a \rightarrow \infty$, we have
\beqnn
\limsup_{n\rightarrow\infty}g(n)
\le \limsup_{u\to\infty}g(u)\lim_{n\rightarrow\infty}\mbf{P}_n\{\kappa_\infty^+<\infty\}.
\eeqnn
It follows that
\[\lim_{n\rightarrow\infty}\mbf{P}_n\{\kappa_\infty^+<\infty\}
=1. \]
Then for large enough initial value $n_0$, the process $Z$ has a positive probability to explode. The following is inspired by \cite[Proof of Theorem 3.1]{FZ22}. Since $Z$ is irreducible in $\mbb{N}_+$, its probability of explosion starting from $1$ is also positive, i.e. $\mbf{P}_1\{\kappa_\infty^+ < \infty\} > 0.$ We take $t > 0$ such that $\mbf{P}_1\{\kappa_\infty^+ \le t\} > 0.$ The stochastic monotonicity in the initial states, see Lemma \ref{Fl3.4}, ensures that for any $n \ge 1$, $\mbf{P}_{n}\{\kappa_\infty^+ > t\} \le \mbf{P}_1\{\kappa_\infty^+ > t\} < 1.$ Let $m \ge 2$, by the Markov property at time $(m - 1)t$, we have
\beqnn
\mbf{P}_{n}\{\kappa_\infty^+ > mt\} \ar=\ar \mbf{P}_{n}\{\kappa_\infty^+ > (m-1)t\}\mbf{E}\left[\mbf{P}_{Z_{(m-1)t}^{(n)}}\{\kappa_\infty^+ > t\}\right]\crr
\ar\le\ar  \mbf{P}_{n}\{\kappa_\infty^+ > (m-1)t\} \mbf{P}_1\{\kappa_\infty^+ > t\}.
\eeqnn
By induction,
\beqnn
\mbf{P}_{n}\{\kappa_\infty^+ > mt\} \le [\mbf{P}_1\{\kappa_\infty^+ > t\}]^m \rightarrow 0
\eeqnn
as $m \rightarrow \infty$. Therefore, $\mbf{P}_{n}\{\kappa_\infty^+ < \infty\} = 1$ for any $n \ge 1$.
Process $Z$ thus explodes.
\qed

By a modification of \cite[Lemma 3.2]{MYZ21}, we have the following result.

\begin{proposition}\label{stay_infinite}
Suppose that there exist a bounded  strictly positive function $g$ on $\mathbb{N}_+$ and an eventually strictly positive function $d$ on $\mathbb{N}_+$ such that $\lim_{n\rightarrow\infty}g(n)=0$ and for large $a > 0$,
\beqlb\label{eq_0412}
\mcr{L}g(n)\leq d(a)g(n)
\eeqlb
holds for all $n>a$. Then
for all $a, t >0$,
\beqlb\label{2.4}
\lim_{n\to\infty}\mbf{P}_n\{\kappa^-_a<t\wedge\kappa_\infty^+\}=0.
\eeqlb
\end{proposition}

\proof
For any $0 <a<n<b$, by \eqref{eq_0412} we have
\beqnn
\mbf{E}_n\big[g(Z_{t\wedge \kappa_{a,b}})\big] \ar=\ar g(n)
+\mbf{E}_n\left[\int_0^{t\wedge \kappa_{a,b}}\mcr{L}g(Z_s)\dd s\right]\cr
\ar\le\ar g(n)
+d(a)\mbf{E}_n\left[\int_0^{t\wedge \kappa_{a,b}-}g(Z_s)\dd s\right]       \cr
\ar\le\ar
g(n)+d(a)\int_0^t\mbf{E}_n\big[g(Z_{s\wedge \kappa_{a,b}})\big]\dd s.
\eeqnn
It follows from Gronwall's lemma that
\beqlb\label{eq0318}
\mbf{E}_n\big[g(Z_{t\wedge \kappa_{a,b}})\big]
\le
g(n)\e^{d(a)t}.
\eeqlb
Then by Fatou's lemma,
\beqnn
\mbf{E}_n\big[g(Z_{\kappa^-_a})1_{\{\kappa^-_a<t\wedge\kappa_\infty^+\}}\big]
\ar\le\ar
\liminf_{b\to\infty}
\mbf{E}_n\big[g(Z_{\kappa^-_a})1_{\{\kappa^-_a<t\wedge\kappa_b^+\}}\big]\cr
\ar\le\ar
\liminf_{b\to\infty}\mbf{E}_n\big[g(Z_{t\wedge \kappa_{a,b}})\big]
\le g(n)\e^{d(a)t},
\eeqnn
which implies
\beqnn
\inf\{g(y): 1\leq y\leq a\}\mbf{P}_n\{\kappa^-_a<t\wedge\kappa_\infty^+\}
\ar\le\ar \mbf{E}_n[g(Z_{\kappa^-_a})1_{\{\kappa^-_a<t\wedge\kappa_\infty^+\}}] \le g(n)\e^{d(a)t}.
\eeqnn
Recall $g$ is a bounded and eventually strictly positive function on $\mbb{N}_+$ with $\lim_{n\to\infty}g(n)=0$. Letting $n \rightarrow \infty$ in above, the result follows.
\qed

\subsection{Criteria for block counting processes on $\bar{\mbb{N}}_+$}  Recall that $N$ is the block counting process of $\Pi$, which is a Markov process on $\bar{\mbb{N}}_+$, see \cite[Theorem 2.3]{FZ23}. In this subsection, we present several criteria that determine the behavior of the process $N$, including conditions for non-explosion, explosion, staying infinite and coming down from infinity. The results in Subsection \ref{s3} play a crucial role in the following proof.

The process $(N^{(n)}_t, t < \tau_\infty^{n,+})$, introduced in Subsection \ref{s3}, is an irreducible Markov process on $\mbb{N}_+$ with initial value $n$ and generator $\mathcal{L}$ given by \eqref{L}, where $\tau_\infty^{n,+} = \inf\{t > 0: N_{t-}^{(n)} = \infty\}$. For any positive constants $a, b$ with $a < n < b$, recall that $\tau_a^{n,-} := \inf\{t > 0: N_t^{(n)} \le a\}$, let $\tau_b^{n,+} := \inf\{t > 0: N_t^{(n)} \ge b\}$ and $\tau_{a,b}^n := \tau_a^{n,-}\wedge\tau_b^{n,+}$. For simplicity, we denote them as $\tau_\infty^+$, $\tau_a^-$, $\tau_b^+$ and $\tau_{a,b}$ when the initial value is specified. It is easy to see that the results in Subsection \ref{ss4.1} hold for $(N^{(n)}_t, t < \tau_\infty^{n,+})$.

Notice that the process $N$ does not explode under the assumption of Proposition \ref{p0607}. Conversely, we establish the following complementary result.

\begin{proposition}\label{c0407}
Under the assumptions of Propositions \ref{CDI} and \ref{stay_infinite}, the process $N$ explodes.
\end{proposition}

\proof
Notice that $\lim_{n\to\infty}\mbf{P}_n\{\tau^-_a<t\wedge\tau_\infty^+\}=0$ for all $a, t > 0$ by Proposition \ref{stay_infinite}. Then the result follows from the second item of Proposition \ref{CDI}.
\qed

\begin{proposition}\label{stay_infinite1}
Under the assumption of Proposition \ref{stay_infinite}, the process $N$ stays infinite.
\end{proposition}

\proof If the process $N$ does not explode, then $\tau^{n,+}_\infty=\infty$ almost surely. By Proposition \ref{stay_infinite},
\beqnn
\lim_{n\to\infty}\mbf{P}_n\{\tau^-_a<t\}
=
\lim_{n\to\infty}\mbf{P}_n\{\tau^-_a<t\wedge\tau_\infty^+\}=0
\eeqnn
for any $t, a >0$.
Moreover, by Lemma \ref{Fl3.4}, the stochastic monotonicity in the initial states, for any $t, a > 0$ we have
\beqnn
\mbf{P}_\infty\{\tau_a^- < t\} \le \mbf{P}_n\{\tau_a^- < t\} \rightarrow 0 \quad \text{as}\ n \rightarrow \infty.
\eeqnn
Then the process $N$ stays infinite.

For the case of explosion, the proof is inspired by \cite[Lemma 3.18]{F22}. Assume that $N$ comes down from infinity. There exists a contradiction with the assumption. In fact, by the Zero-One law in \cite[Lemma 2.5]{F22}, the process leaves infinity instantaneously. Recall that $\tau_\infty^{\infty, +} = \inf\{t > 0: N_{t-}^{(\infty)} = \infty\}$. We consider an excursion from $\infty$ with length $\tau_\infty^{\infty, +}$, such that $\tau_\infty^{\infty, +} > \tau_a^{\infty,-} > t$ for some $t > 0$ and $a > 0.$ By the Markov property at time $t$, conditionally on $N_t$, the process $(N_{t + s}^{(\infty)}, 0 \le s \le \tau_\infty^{\infty, +} - t)$ has the same law as the process $(N_s^{(N_t)}, 0 \le s \le \tau_\infty^{N_t,+})$, where $\tau_\infty^{N_t,+} := \inf\{s > 0: N_{s-}^{(N_t)} = \infty\}.$ We extend  the definition of $g$ at $\infty$ as $g(\infty) = \lim_{n \rightarrow \infty}g(n) = 0$. Then $g$ is a bounded and positive function on $\bar{\mbb{N}}_+$ with $g(\infty) = 0$ and $g(n) > 0$ for any $n \in \mbb{N}_+$. Similar to \eqref{eq0318}, by the assumption of Proposition \ref{stay_infinite}, we have
\beqnn
\mbf{E}_\infty\left[g(N_{(t + s)\wedge \tau_{a}^-\wedge \tau_\infty^+})1_{\{t+s < \tau_a^- < \tau_\infty^+\}}|N_t\right] \ar=\ar \mbf{E}_{N_t}\left[g(N_{s\wedge \tau_{a}^-\wedge \tau_\infty^{+}})1_{\{s < \tau_a^- < \tau_\infty^{+}\}}\right]\cr
\ar\le\ar \mbf{E}_{N_t}\big[g(N_{s\wedge \tau_{a}^-\wedge \tau_\infty^{+}})\big]\cr
\ar\le\ar g(N_t)\e^{d(a)s}.
\eeqnn
Letting $t \rightarrow 0+$, we have
\beqnn
\lim_{t \rightarrow 0+}\mbf{E}_\infty\left[g(N_{(t + s)\wedge \tau_{a}^-\wedge \tau_\infty^+})1_{\{t+s < \tau_a^- < \tau_\infty^+\}}\big|N_t\right] = \mbf{E}_\infty\left[g(N_{s\wedge \tau_{a}^-\wedge \tau_\infty^+})1_{\{s < \tau_a^- < \tau_\infty^+\}}\right] = 0.
\eeqnn
Then for any $s \in (0,  \tau_{a}^{\infty,-})$, $N_s^{(\infty)} = \infty$ almost surely. This contradicts the fact that the process leaves infinity instantaneously. As a result, the process cannot come down from infinity, i.e., the process stays infinite.
\qed

Recall that $(N_m^{(n)}(t), t \ge 0)$ is a Markov process on $\mbb{N}_+$ with initial value $N_m^{(n)}(0) = n$ and generator $\mathcal{L}^m$ given by \eqref{Lm}. We then have the following result.

\begin{proposition}\label{l0406}
Suppose that
\beqlb\label{eq_0509a}
\lim_{a \rightarrow \infty}\lim_{n\rightarrow\infty}\mbf{P}_n\{\tau_{a,m}^- < t\} = 1
\eeqlb
for any $t > 0$ and $m \in \mathbb{N}_+.$ Then the process $N$ comes down from infinity.
\end{proposition}

\proof
Notice that $N_m$ does not explode and $\lim_{n \rightarrow \infty}\tau_{a,m}^{n,-} = \tau_{a,m}^{\infty,-}$ by Lemmas \ref{Fl3.8} and \ref{Fl3.10}. Then by \eqref{eq_0509a} and dominated convergence theorem, we have
\beqlb\label{eq_0509b}
\lim_{a \rightarrow \infty}\mbf{P}_\infty\{\tau_{a,m}^- < t \} = 1
\eeqlb
for any $t > 0$ and $m \in \mbb{N}_+$, which implies that $\infty$ is an entrance boundary for the process $N_m.$ By the Zero-One law in \cite[Lemma 2.5]{F22}, $N_m$ comes down from infinity immediately and does not explode, i.e., for any $m \in \mathbb{N}$ and $t > 0,$ $N_m^{(\infty)}(t) < \infty$ almost surely. Again by Lemma \ref{Fl3.10}, we have  $\lim_{m \rightarrow \infty}\tau_{a,m}^{\infty,-} = \tau_{a}^{\infty,-}$ almost surely. It follows from \eqref{eq_0509b} and dominated convergence theorem that
$$\lim_{a \rightarrow \infty}\mbf{P}_\infty\{\tau_a^- < t \} = 1$$
for any $t > 0.$  Then there exists a constant $a$ big enough such that $\mbf{P}_\infty\{\tau_a^- < t \} > 0$ for any $t > 0.$ We take a proper $t$ such that $\mbf{P}_\infty\{\tau_a^- > t \} \in (0, 1)$.  Then for $n \ge 2$  combining the Markov property at $(n-1)t$,  Lemma \ref{Fl3.4}, and the stochastic monotonicity in the initial states, we have
\beqnn
\mbf{P}_\infty\{\tau_a^- > nt\} \ar=\ar \mbf{P}_\infty\{\tau_a^- > (n-1)t\} \mbf{E}\left[\mbf{P}_{N_{(n-1)t}^{(\infty)}}\{\tau_a^- > t\}\right]\crr
\ar\le\ar \mbf{P}_\infty\{\tau_a^- > (n-1)t\}\mbf{P}_\infty\{\tau_a^- > t\}.
\eeqnn
By induction, one sees that
\beqnn
\mbf{P}_\infty\{\tau_a^- > nt\} \le \left[\mbf{P}_\infty\{\tau_a^- > t\} \right]^n \rightarrow 0
\eeqnn
as $n \rightarrow \infty$. Then $\mbf{P}_\infty\{\tau_a^- < \infty \} = 1$, i.e., the process $N$ comes down from infinity. The result follows.
\qed

\section{Proofs of the main results}

{\bf Proof of Theorem \ref{p3}.}
It suffices to show that (\ref{F-L}) holds. For $m \in \mbb{N}_+$, set $g(n) = \log^{(m)} (n+l)$, where $l$ is a constant such that $\log^{(m-1)} l > 0$. By the inequality  $\log (1-x) \le -x$ for any $x \in (0, 1)$, we have
\beqlb\label{eq0519}
\ar\ar \log\log(n - k + l + 1) - \log\log (n + l) \cr
\ar\ar\quad = \log\left(\frac{\log (n - k + l + 1)}{\log (n + l)}\right)\cr
\ar\ar\quad  = \log\left(1 + \frac{\log(1-(k-1)/(n+l))}{\log (n + l) }\right)\cr
\ar\ar\quad  \le\frac{\log(1-(k-1)/(n + l))}{\log (n + l)}\cr
\ar\ar\quad  \le -\frac{k - 1}{(n + l) \log (n + l)}.
\eeqlb
Moreover,
\beqnn
\ar\ar \log^{(m)}(n - k + l + 1) - \log^{(m)} (n+l)\crr
\ar\ar\quad = \log\left(1 + \frac{\log^{(m-1)}(n - k + l+1) - \log^{(m-1)} (n + l)}{\log^{(m-1)} (n + l)}\right)\crr
\ar\ar\quad \le \frac{\log^{(m-1)}(n - k + l + 1) - \log^{(m-1)} (n + l)}{\log^{(m-1)} (n + l)}.
\eeqnn
Then by induction, one obtains that
\beqnn
\log^{(m)}(n - k + l + 1) - \log^{(m)} (n + l) \le -\frac{k-1}{\prod_{\ell = 0}^{m - 1} \log^{(\ell)} (n+l)}.
\eeqnn
Here $\prod_{\ell = 1}^0 := 1.$ By the above and Proposition \ref{p0217}, as $n \rightarrow \infty$ we have
\beqlb\label{Lc_slow}
\mathcal{L}^c\log^{(m)} (n + l)
\ar=\ar\sum_{k=2}^n{n\choose k}\lambda_{n,k}[\log^{(m)}(n-k+ l + 1)-\log^{(m)} (n + l)]\cr
\ar\le\ar -\frac{1}{\prod_{\ell = 0}^{m - 1} \log^{(\ell)}(n + l)}\sum_{k=2}^n{n\choose k}\lambda_{n,k} (k-1)\cr
\ar=\ar -\frac{\Phi_\Lambda(n)}{\prod_{\ell = 0}^{m - 1} \log^{(\ell)}(n + l)}.
\eeqlb

On the other hand, by inequality $\log (1+x)\le x$ for any $x > 0$, we have
\beqnn
\ar\ar\log\log (n + l + k) - \log\log (n + l)\cr 
\ar\ar\quad = \log\left(1 + \frac{\log (1 + k/(n + l))}{\log (n + l)}\right)\cr
\ar\ar\quad \le \frac{\log (1 + k/(n + l))}{\log (n + l)}
\eeqnn 
and
\beqnn
\ar\ar \log^{(m)} (n + l + k) - \log^{m} (n + l) \cr 
\ar\ar\quad = \log\left(1 + \frac{\log^{(m-1)} (n + l + k) - \log^{(m - 1)} (n + l)}{\log^{(m - 1)} (n + l)}\right)\crr
\ar\ar\quad \le \frac{\log^{(m-1)} (n + l + k) - \log^{(m - 1)} (n + l)}{\log^{(m - 1)} (n + l)}.
\eeqnn 
It follows from induction that
\beqlb\label{eq0417}
\log^{(m)} (n + l + k) - \log^{m} (n + l) \le \frac{\log (1 + k/(n + l))}{\prod_{\ell = 1}^{m - 1} \log^{(\ell)} (n + l)}.
\eeqlb
Moreover, similar to the proofs of  Propositions \ref{p0211}, \ref{mu_tail}, one sees that
\beqnn
n\sum_{k=1}^\infty\mu(k) \log \left(1 + \frac{k}{n + l}\right) = \frac{\Phi_\mu(n)}{n + l} + O((\log n)^\alpha)
\eeqnn
as $n \rightarrow \infty$. Then by the above, \eqref{eq0417} and Condition \ref{C1}, we have
\beqlb\label{Lf_slow}
\mathcal{L}^f\log^{(m)} (n + l) \ar=\ar n\sum_{k=1}^\infty\mu(k)\left[\log^{(m)}(n+l+k)-\log^{(m)} (n+l)\right] \cr
\ar\le\ar \frac{n}{\prod_{\ell = 1}^{m - 1} \log^{(\ell)} (n+l)} \sum_{k=1}^\infty\mu(k) \log \left(1 + \frac{k}{n+l}\right) \cr
\ar =\ar \frac{\Phi_\mu(n)}{\prod_{\ell = 0}^{m - 1} \log^{(\ell)} (n + l)} + O\left(\frac{(\log n)^{\alpha-1}}{\prod_{\ell = 2}^{m - 1} \log^{(\ell)} n}\right)
\eeqlb
as $n \rightarrow \infty.$
It follows from \eqref{eq0120} that there exists a constant $C > 0$ such that for large $n$,
\beqnn
\frac{ \Phi_\Lambda(n)-\Phi_\mu(n)}{\prod_{\ell = 0}^m \log^{(\ell)} (n + l) } \ge -C
\eeqnn
holds.  Recall that $\alpha \in (0, 1]$. 
By the above, \eqref{Lc_slow} and \eqref{Lf_slow}, we have
\beqnn
\mathcal{L}\log^{(m)} (n + l)\ar=\ar\mathcal{L}^c\log^{(m)} (n + l)+\mathcal{L}^f\log^{(m)} (n + l)\cr
\ar\le\ar	
- \frac{\Phi_\Lambda(n)- \Phi_\mu(n)}{\prod_{\ell = 0}^{m - 1} \log^{(\ell)} (n + l)} +O\left(\frac{(\log n)^{\alpha-1}}{\prod_{\ell = 2}^{m - 1} \log^{(\ell)} n}\right)\cr
\ar\le\ar C \log^{(m)} (n + l)  + O\left(\frac{(\log n)^{\alpha-1}}{\prod_{\ell = 2}^{m - 1} \log^{(\ell)} n}\right)
\eeqnn 
as $n \rightarrow \infty$.
Then there exists a constant $C > 0$ such that \eqref{F-L} holds for all $n \ge 1$. The result follows by Proposition \ref{p0607}.
\qed

{\bf Proof of Theorem \ref{pp0406}.}
For a fixed constant $l \ge 10$, we choose $g_l(n) = 1 + \frac{1}{\log\log(n + l)}$. Then $\lim_{n \rightarrow \infty}g_l(n) = 1$ and $\sup_n g_l(n) = 1 + \frac{1}{\log\log(l+1)}.$
Then by \eqref{eq0519} and the fact of $\log(1 + x) \le x$ for any $x \in (-1, 1),$ we have
\beqlb\label{Lc0514}
\mathcal{L}^c\left(1+\frac{1}{\log\log(n +l)}\right)
\ar=\ar \sum_{k=2}^n{n\choose k}\lambda_{n,k}[\log\log(n-k+1+l)^{-1} - \log\log(n+l)^{-1}]\cr
\ar=\ar -\sum_{k=2}^n{n\choose k}\lambda_{n,k}
\frac{\log\left(1 + \frac{\log(1-(k-1)/(n+l))}{\log(n+l)}\right) }{\log\log(n+l)\log\log(n-k+1+l)}\cr
\ar\ge\ar  -\sum_{k=2}^n{n\choose k}\lambda_{n,k}
\frac{\log(1-\frac{k-1}{n+l})}{\log(n+l)\log\log(n+l)\log\log(n-k+1+l)}\cr
\ar\geq\ar -\frac{1}{\log(n+l)(\log\log(n+l))^2}\sum_{k=2}^n{n\choose k}\lambda_{n,k} \log\left(1-\frac{k-1}{n+l}\right)\cr
\ar\geq\ar \frac{1}{(n+l)\log(n+l)(\log\log (n+l))^2}\sum_{k=2}^n{n\choose k}\lambda_{n,k}(k-1)\cr
\ar=\ar \frac{\Phi_\Lambda(n)}{n\log n(\log\log n)^2} +O(1)
\eeqlb
as $n \rightarrow \infty,$
where $O(1)$ is due to changing $n+l$ to $n$. On the other hand, by the fact of $\log(1 + x)\le x$ for any $x \in (0, 1)$, Propositions \ref{p0211} and \ref{mu_tail}, as $n \rightarrow \infty,$ we have
\beqlb\label{Lf}
\mathcal{L}^f\left(1 + \frac{1}{\log\log(n+l)}\right) \ar=\ar n\sum_{k = 1}^\infty \mu(k)[(\log\log(n + l+ k))^{-1} - (\log\log(n+l))^{-1}]\cr
\ar=\ar n\sum_{k = 1}^\infty\mu(k)\frac{\log\log(n+l) - \log\log(n +l+ k)}{\log\log(n+l)\log\log(n+l+k)}\cr
\ar\ge\ar -\frac{n}{(\log\log n)^2}\sum_{k = 1}^\infty \mu(k)[\log\log(n + l+ k) - \log\log(n+l)]\cr
\ar=\ar -\frac{n}{(\log\log n)^2}\sum_{k = 1}^\infty \mu(k)\log\left(1 + \frac{\log(1 + k/(n+l))}{\log(n+l)}\right)\cr
\ar\ge\ar -\frac{n}{(\log\log n)^2}\sum_{k = 1}^\infty \mu(k)\frac{\log(1 + k/(n+l))}{\log(n+l)}\cr
\ar\ge\ar  -\frac{n}{\log n(\log\log n)^2}\sum_{k = 1}^\infty \mu(k)\log\left(1 + \frac{k}{n}\right)\cr
\ar=\ar -\frac{\Phi_\mu(n)}{n\log n(\log\log n)^2} -\left(2\log2 - \frac{1}{2}\right)\frac{b(\log n)^{\alpha - 1}}{(\log\log n)^2} \cr
\ar\ar +  o\left(\frac{(\log n)^{\alpha -1}}{(\log\log n)^2}\right).
\eeqlb
Recall that $\alpha \in (0, 1].$ By  Lemma \ref{ll0406}, \eqref{Lc0514} and \eqref{Lf}, for any $m \ge 1$ one sees that
\beqlb\label{eq0208a}
\mathcal{L}^m\left(1 + \frac{1}{\log\log(n + l)}\right) \ar\ge\ar \mathcal{L}\left(1 + \frac{1}{\log\log(n + l)}\right) \cr
\ar=\ar \mathcal{L}^c\left(1 + \frac{1}{\log\log(n + l)}\right) + \mathcal{L}^f\left(1 + \frac{1}{\log\log(n + l)}\right)\cr
\ar\ge\ar \frac{\Phi_\Lambda(n) - \Phi_\mu(n)}{n\log n (\log\log n )^2} - \left(2\log2 - \frac{1}{2}\right)\frac{b(\log n)^{\alpha - 1}}{(\log\log n)^2} \cr
\ar\ar + o\left(\frac{(\log n )^{\alpha - 1}}{(\log\log n )^2}\right)
\eeqlb
as $n \rightarrow \infty$,
which goes to infinity by \eqref{Phi_0518}. Then there exists an eventually strictly positive function $d(a)$ on $\mbb{N}_+$ with $\lim_{a \rightarrow \infty}d(a) = \infty$ such that \eqref{L0208} holds. Moreover, the process $N_m$ does not explode by Lemma \ref{Fl3.8}. Then by Proposition \ref{CDI}, we obtain that
\beqnn
\lim_{a \rightarrow \infty}\lim_{n \rightarrow \infty}\mbf{P}_n\{\tau_{a,m}^- < t\}\ge \frac{\limsup_{u \rightarrow \infty}g_l(u)}{\sup_u g_l(u)} = \frac{\log\log(l+1)}{\log\log(l + 1) + 1}.
\eeqnn
Letting $l \rightarrow \infty$, for any $t > 0$ and $m \in \mathbb{N}$, we have
\beqnn
\lim_{a \rightarrow \infty}\lim_{n \rightarrow \infty}\mbf{P}_n\{\tau_{a,m}^- < t\} = 1.
\eeqnn
The result follows from Proposition \ref{l0406}.
\qed

\begin{lemma}\label{L^c4}
If Condition \ref{C2} holds, then as $n\to\infty$,
\beqnn
\mathcal{L}^c\left(\frac{1}{\log\log (n+10)}\right) \ar\le\ar \frac{\Phi_\Lambda(n)}{(n+10)\log (n+10) (\log\log (n + 10))^2} \cr
\ar\ar +  \frac{2C_2}{\beta  } + O\left(\frac{1}{\log\log n}\right),
\eeqnn
where $C_2$ is the constant given in \eqref{Lambda}.
\end{lemma}

\proof
For any $2 \le k \le n$, by \eqref{eq0519} we have
\beqlb\label{eq0519a}
\frac{\log\log (n+10)}{\log\log(n - k + 11)} \ar=\ar 1 - \frac{\log\log(n - k + 11) - \log\log(n + 10)}{\log\log(n - k + 11)}\cr
\ar=\ar 1 - \frac{\log\left(1 + \frac{\log(1 - (k-1)/(n+10))}{\log (n+10)}\right)}{\log\log(n - k + 11)}\cr
\ar\le\ar 1 - \frac{\log\left(1 + \frac{\log(1 - (k-1)/(n+10))}{\log (n+10)}\right)}{\log\log(11)}.
\eeqlb
Then by the definition of $\mathcal{L}^c$ and \eqref{eq0519a}, one sees that
\beqlb\label{eq0516}
\mathcal{L}^c\left(\frac{1}{\log\log(n+10)}\right)  
\ar  =\ar \sum_{k=2}^n{n\choose k}\lambda_{n,k}\left[\frac{1}{\log\log(n-k+11)}-\frac{1}{\log\log(n+10)}\right]\cr
\ar  =\ar \sum_{k=2}^n{n\choose k}\lambda_{n,k}
\frac{\log\log(n+10)-\log\log(n-k+11)  }{\log\log(n+10)\log\log(n-k+11)}\cr
\ar  =\ar  -\sum_{k=2}^n{n\choose k}\lambda_{n,k}
\frac{\log\left(1 + \frac{\log(1-(k-1)/(n+10))}{\log (n+10)}\right)  }{(\log\log(n+10))^2}\cdot\frac{\log\log (n+10)}{\log\log(n - k + 11)}\cr
\ar \le\ar -\frac{1}{(\log\log (n+10))^2}\sum_{k = 2}^n{n\choose k}\lambda_{n,k}
\log\left(1 + \frac{\log(1-(k-1)/(n+10))}{\log (n+10)}\right)\cr
\ar\ar\quad  + \frac{1}{(\log\log (n+10))^2}\sum_{k = 2}^n{n\choose k}\lambda_{n,k} \frac{\left(\log\left(1 + \frac{\log(1 - (k-1)/(n+10))}{\log (n+10)}\right)\right)^2}{\log\log(11)} \cr
\ar  =: \ar  I_1(n) + I_2(n),
\eeqlb
where the inequality above comes from \eqref{eq0519a}.

Now we recall some elementary inequalities about $\log(1 - x)$. For any fixed $\delta \in (0, 1)$, there exist constants $C_\delta, \tilde{C}_\delta > 0$ depending on $\delta$ such that
\beqlb\label{log_eq1}
\log(1 - x) \ge - x - C_\delta x^2, \quad (\log(1 - x))^2 \le \tilde{C}_\delta x^2, \qquad \forall\ x \in (0, \delta).
\eeqlb
Fix a constant $\delta \in (0,1)$ satisfying
\beqlb\label{delta}
\frac{\delta^{-\beta}}{\log\log (11)} \le 2,
\eeqlb
which is feasible since $2\log\log(11) > 1,$
it follows from \eqref{log_eq1} and Proposition \ref{Phi2} that, as $n\to\infty$,
\beqlb\label{I_11}
\ar\ar -\frac{1}{(\log\log (n+10))^2}\sum_{k = 2}^{\lfloor\delta n\rfloor}{n\choose k}\lambda_{n,k}
\log\left(1 + \frac{\log(1-(k-1)/(n+10))}{\log (n+10)}\right)\cr
\ar\ar\quad\le -\frac{1}{\log (n+10) (\log\log (n+10))^2}\sum_{k = 2}^{\lfloor\delta n\rfloor} {n\choose k}\lambda_{n,k} \log\left(1 - \frac{k - 1}{n+10}\right)\cr
\ar\ar \qquad+ \frac{C_\delta}{(\log (n+10) \log\log (n+10))^2}\sum_{k = 2}^{\lfloor\delta n\rfloor} {n\choose k}\lambda_{n,k} \left(\log\left(1 - \frac{k - 1}{n+10}\right)\right)^2\cr
\ar\ar\quad\le -\frac{1}{\log (n+10) (\log\log (n+10))^2}\sum_{k = 2}^{n} {n\choose k}\lambda_{n,k} \log\left(1 - \frac{k - 1}{n+10}\right)\cr
\ar\ar \qquad+ \frac{C_\delta}{(\log (n+10) \log\log (n+10))^2}\sum_{k = 2}^{\lfloor\delta n\rfloor} {n\choose k}\lambda_{n,k} \left(\log\left(1 - \frac{k - 1}{n+10}\right)\right)^2\cr
\ar\ar\quad\le \frac{\Phi_\Lambda(n)}{(n+10)\log (n+10) (\log\log (n+10))^2} + O\left(\frac{1}{\log n (\log\log n)^2}\right)\cr
\ar\ar\qquad + \frac{C_\delta \cdot \tilde{C}_\delta}{((n+10)\log (n+10) \log\log (n+10))^2} \sum_{k = 2}^{n} {n\choose k}\lambda_{n,k} \left(k - 1\right)^2\cr
\ar\ar\quad= \frac{\Phi_\Lambda(n)}{(n + 10)\log (n+10) (\log\log (n + 10))^2} + O\left(\frac{1}{\log n (\log\log n)^2}\right).
\eeqlb
where the first term in the third inequality arises from
\beqnn
\sum_{k=2}^n{n\choose k}\lambda_{n,k} \log \left(1 - \frac{k - 1}{n+10}\right) = -\frac{\Phi_\Lambda(n)}{n+10} + O(1) 
\eeqnn
as $n\to\infty$, which can be obtained similarly to the proof of Proposition \ref{p0217}.
Moreover, by Proposition \ref{l0225}, as $n \rightarrow \infty,$ we have
\beqlb\label{I_12}
\ar\ar -\frac{1}{(\log\log (n+10))^2}\sum_{k = \lfloor\delta n\rfloor + 1}^{n}{n\choose k}\lambda_{n,k}
\log\left(1 + \frac{\log(1-(k-1)/(n+10))}{\log (n+10)}\right)\cr
\ar\ar\quad \le \frac{1}{\log\log (n + 10)}\sum_{k = \lfloor\delta n\rfloor + 1}^{n}{n\choose k}\lambda_{n,k} \le \frac{C_2 \delta^{-\beta}/\beta}{\log\log n} +  O\left(\frac{1}{n\log\log n}\right).
\eeqlb
Then by \eqref{eq0516}, \eqref{I_11} and \eqref{I_12}, one sees that, as $n\to\infty$,
\beqlb\label{I_1}
I_1(n) \le \frac{\Phi_\Lambda(n)}{(n+10)\log (n+10) (\log\log (n + 10))^2} + O\left(\frac{1}{\log\log n}\right).
\eeqlb
On the other hand, by \eqref{log_eq1} and Proposition \ref{Phi2},
\beqlb\label{I21}
\ar\ar\frac{1}{(\log\log (n+10))^2}\sum_{k = 2}^{\lfloor\delta n\rfloor}{n\choose k}\lambda_{n,k} \left(\log\left(1 + \frac{\log(1 - (k-1)/(n+10))}{\log (n+10)}\right)\right)^2 \cr
\ar\ar\quad\le \frac{\tilde{C}_\delta}{ (\log(n+10)\log\log (n+10))^2}\sum_{k = 2}^{\lfloor\delta n\rfloor}{n\choose k}\lambda_{n,k}\left( \log\left(1 - \frac{k-1}{n+10}\right)\right)^2\cr
\ar\ar\quad\le \frac{\tilde{C}_\delta^2}{(n\log n \log\log n)^2}\sum_{k = 2}^n {n\choose k}\lambda_{n,k}(k - 1)^2 \cr
\ar\ar\quad=  O\left(\frac{1}{(\log n\log\log n)^2}\right)
\eeqlb
as $n \rightarrow \infty.$
Moreover, by Proposition \ref{l0225}, as $n \rightarrow \infty$, one obtains that
\beqlb\label{I22}
\ar\ar \frac{1}{(\log\log (n+10))^2}\sum_{k = \lfloor\delta n\rfloor + 1}^{n}{n\choose k}\lambda_{n,k} \left(\log\left(1 + \frac{\log(1 - (k-1)/(n+10))}{\log (n+10)}\right)\right)^2 \cr
\ar\ar\qquad\le \sum_{k = \lfloor\delta n\rfloor + 1}^{n}{n\choose k}\lambda_{n,k} \le \frac{C_2}{\beta} \delta^{-\beta} + O(n^{-1}).
\eeqlb
Then by \eqref{eq0516}, \eqref{I21}, \eqref{I22} and \eqref{delta}, as $n \rightarrow \infty$ we have
\beqlb\label{I_2}
I_2(n) \ar\le\ar \frac{C_2 \delta^{-\beta}}{\beta\log\log (11)} +  O\left(\frac{1}{(\log n\log\log n)^2}\right)\cr
\ar\le\ar \frac{2C_2}{\beta} + O\left(\frac{1}{(\log n\log\log n)^2}\right).
\eeqlb
Then the result follows from  \eqref{eq0516}, \eqref{I_1} and \eqref{I_2}.
\qed

\begin{lemma}\label{L^f4}
Assume that Condition \ref{C1} holds. Then as $n\to\infty$, we have
\beqnn
\mathcal{L}^f \left(\frac{1}{\log\log (n+10)}\right)\ar\le\ar -\frac{ \Phi_\mu(n)}{(n+10)\log (n+10)(\log\log (n+10))^2} \cr
\ar\ar  - \left(2\log 2 - \frac{1}{2}\right)\frac{b(\log n)^{\alpha - 1}}{(\log\log n)^2} + o\left(\frac{(\log n)^{\alpha - 1}}{(\log\log n)^2}\right).
\eeqnn
\end{lemma}

\proof
Notice that
\beqlb\label{Lf0513}
\mathcal{L}^f \left(\frac{1}{\log\log (n+10)}\right)\ar=\ar n\sum_{k = 1}^\infty \mu(k)\left[\frac{1}{\log\log(n+k+10)} - \frac{1}{\log\log (n+10)}\right] \cr
\ar=\ar - n\sum_{k = 1}^n \mu(k)\frac{\log\left(1 + \frac{\log(1 + k/(n+10))}{\log (n+10)}\right)}{\log\log(n+k+10)\log\log(n+10)}\cr
\ar\ar - n\sum_{k = n+1}^\infty \mu(k)\frac{\log\left(1 + \frac{\log(1 + k/(n+10))}{\log (n+10)}\right)}{\log\log(n+k+10)\log\log(n+10)} \cr
\ar=:\ar \tilde{I}_1(n) + \tilde{I}_2(n).
\eeqlb
By the fact of
\beqnn
x \ge\log(1 + x) \ge x - \frac{1}{2}x^2, \qquad \forall\ x \in (0, 1),
\eeqnn
one obtains that
\beqnn
\tilde{I}_1(n)
\ar=\ar - n\sum_{k = 1}^n \mu(k)\frac{\log\left(1 + \frac{\log(1 + k/(n+10))}{\log (n+10)}\right)}{(\log\log(n+10))^2}\cdot \frac{\log\log (n+10)}{\log\log(n+k+10)}\cr
\ar\le\ar -\frac{n}{\log (n+10)(\log\log (n+10))^2} \sum_{k = 1}^n \mu(k)\log\left(1 + \frac{k}{n+10}\right)\cdot\frac{\log\log (n+10)}{\log\log (n+k+10)}\cr
\ar\ar + \frac{n}{2(\log (n+10)\log\log (n+10))^2}\sum_{k = 1}^n \mu(k)\left(\log\left(1 + \frac{k}{n+10}\right)\right)^2.
\eeqnn
Similar to \eqref{eq0519}, we have
\beqlb\label{eq0519b}
\frac{\log\log (n+10)}{\log\log (n+k+10)} \ar=\ar 1- \frac{\log\log (n+k+10) - \log\log (n+10)}{\log\log (n+k+10)}\cr
\ar=\ar 1 - \frac{\log\left(1 + \frac{\log(1 + k/(n+10))}{\log (n+10)}\right)}{\log\log (n + k+10)}\cr
\ar\ge\ar 1 - \frac{\log(1 + k/(n+10))}{\log (n+10)\log\log(n+10)}
\eeqlb
for any $1 \le k \le n.$ Recall the fact of $(\log (1 + x))^2 \le x^2$
for any $x \in (0, 1)$.
By the proof of Proposition \ref{p0211}, as $n \rightarrow \infty$, one sees that
\beqnn
n\sum_{k = 1}^n \mu(k)\log\left(1 + \frac{k}{n+10}\right) = \frac{\Phi_\mu(n)}{n+10} - \frac{b(\log n)^\alpha}{2} + o((\log n)^\alpha).
\eeqnn
Then by the above,
\eqref{mu2} and \eqref{eq0519b}, as $n \rightarrow \infty$, one obtains that
\beqlb\label{tildeI1}
\tilde{I}_1(n)
\ar\le\ar -\frac{n}{\log (n+10)(\log\log (n+10))^2} \sum_{k = 1}^n \mu(k)\log\left(1 + \frac{k}{n+10}\right) \cr
\ar\ar + \frac{n(2 + \log\log n)}{2(\log n)^2(\log\log n)^3} \sum_{k = 1}^n \mu(k)\left(\log\left(1 + \frac{k}{n}\right)\right)^2\cr
\ar\le\ar -\frac{ \Phi_\mu(n)}{(n+10)\log (n+10)(\log\log (n+10))^2} + \frac{b(\log n)^{\alpha - 1}}{2(\log\log n)^2}\cr
\ar\ar  + o\left(\frac{(\log n)^{\alpha - 1}}{(\log\log n)^2}\right) + \frac{2+\log\log n}{2n(\log n)^2(\log\log n)^3} \sum_{k = 1}^n \mu(k)k^2  \cr
\ar=\ar  -\frac{ \Phi_\mu(n)}{(n+10)\log (n+10)(\log\log (n+10))^2} + \frac{b(\log n)^{\alpha - 1}}{2(\log\log n)^2}\cr
\ar\ar + o\left(\frac{(\log n)^{\alpha - 1}}{(\log\log n)^2}\right).
\eeqlb
Moreover, similar to \eqref{eq0519b}, we have
\beqlb\label{eq0519c}
\frac{\log\log(n+10)}{\log\log(n+k+10)} \ar=\ar 1 - \frac{\log\left(1 + \frac{\log(1 + k/(n+10))}{\log (n+10)}\right)}{\log\log (n + k+10)}\cr
\ar\ge\ar 1 - \frac{\log\left(1 + \frac{\log(1 + k/(n+10))}{\log (n+10)}\right)}{\log\log (n +10)}.
\eeqlb
Applying  inequalities \eqref{eq0519c}, $\log (1 + x) \ge x - x^2/2$ and $(\log (1 + x))^2 \le x^2$ for any $x \in (0, 1)$ together with the proofs of Propositions \ref{mu_tail} and  \ref{c_mu_tail}, as $n \rightarrow \infty,$ we have
\beqlb\label{tt0304}
\tilde{I}_2(n)
\ar=\ar-\frac{n}{(\log\log(n+10))^2}\!\!\sum_{k = n + 1}^\infty\!\mu(k)\log\left(\!1\! +\! \frac{\log(1 + k/(n+10))}{\log(n+10)}\!\right)\!\frac{\log\log(n+10)}{\log\log(n+k+10)}\cr
\ar\le\ar -\frac{n}{(\log\log(n+10))^2}\sum_{k = n + 1}^\infty\mu(k)\log\left(1 + \frac{\log(1 + k/(n+10))}{\log(n+10)}\right)\cr
\ar\ar + \frac{n}{(\log\log(n+10))^3}\sum_{k = n + 1}^\infty\mu(k)\left(\log\left(1 + \frac{\log(1 + k/(n+10))}{\log(n+10)}\right)\right)^2\cr
\ar\le\ar -\frac{n}{\log(n+10)(\log\log(n+10))^2}\sum_{k = n + 1}^\infty\mu(k) \log\left(1 + \frac{k}{n+10}\right)\cr
\ar\ar + \frac{n(2 + \log\log(n+10))}{2(\log(n+10))^2(\log\log(n+10))^3}\sum_{k = n + 1}^\infty\mu(k) \left(\log\left(1 + \frac{k}{n+10}\right)\right)^2\cr
\ar=\ar -\frac{(2\log2) b(\log n)^{\alpha - 1}}{(\log\log n)^2} + o\left(\frac{(\log n)^{\alpha - 1}}{(\log\log n)^2}\right).
\eeqlb
Then the result follows from \eqref{Lf0513}, \eqref{tildeI1} and \eqref{tt0304}.
\qed

{\bf Proof of Theorem \ref{p4}.}
Take $g(n) = (\log\log (n+10))^{-1}.$ Recall that $\alpha \in (0, 1]$. It follows from Lemmas \ref{L^c4} and \ref{L^f4} that
\beqlb\label{L4}
\mathcal{L} \left(\frac{1}{\log\log (n+10)}\right) \ar=\ar \mathcal{L}^c\left(\frac{1}{\log\log (n+10)}\right)+\mathcal{L}^f \left(\frac{1}{\log\log (n+10)}\right)\cr
\ar\le\ar \frac{\Phi_\Lambda(n) - \Phi_\mu(n)}{(n+10)\log (n+10)(\log\log (n+10))^2}\cr
\ar\ar  + \frac{2C_2 }{\beta }  - \left(2\log 2 - \frac{1}{2}\right)\frac{b(\log n)^{\alpha - 1}}{(\log\log n)^2} \cr
\ar\ar + O\left(\frac{1}{\log\log n}\right) + o\left(\frac{(\log n)^{\alpha - 1}}{(\log\log n)^2}\right)\cr
\ar=\ar \frac{\Phi_\Lambda(n) - \Phi_\mu(n)}{(n+10)\log (n+10)(\log\log (n+10))^2}\cr
\ar\ar  + \frac{2C_2 }{\beta } + O\left(\frac{1}{\log\log n}\right)
\eeqlb
as $n \rightarrow \infty.$
By \eqref{Phi_0126}, for large $a > 0,$ there exists a constant $d(a) > 0$ such that
\beqnn
\frac{\Phi_\Lambda(n) -  \Phi_\mu(n)}{(n+10)\log (n+10) \log\log (n+10) } + \frac{2C_2  }{\beta } \log\log (n+10) \le d(a)
\eeqnn
for any $n \ge a$. Then \eqref{eq_0412} holds. The result follows by Proposition \ref{stay_infinite1}.
\qed

{\bf Proof of Theorem \ref{p0206}.}
Recall that $\alpha \in (0,1].$ By \eqref{L4}, one sees that
\beqnn
\mathcal{L} \left(\frac{1}{\log\log (n+10)}\right) \ar=\ar \mathcal{L}^c\left(\frac{1}{\log\log (n+10)}\right)+\mathcal{L}^f \left(\frac{1}{\log\log (n+10)}\right)\cr
\ar\le\ar \frac{\Phi_\Lambda(n) - \Phi_\mu(n)}{(n+10)\log (n+10)(\log\log (n+10))^2}\cr
\ar\ar  + \frac{2C_2  }{\beta } + O\left(\frac{1}{\log\log n}\right)\crr
\ar=\ar  \frac{ \Phi_\Lambda(n) - \Phi_\mu(n)}{n\log n(\log\log n)^2} + O(1)
\eeqnn
as $n \rightarrow \infty$,
which goes to $-\infty$ by \eqref{Phi_0518a}. Then \eqref{eq_0412} holds for large $n$.

Moreover,
\beqnn
\mathcal{L}\left(1 - \frac{1}{\log \log (n + 10)}\right)\ar=\ar -\left[\mathcal{L}^c\left(\frac{1}{\log \log (n + 10)}\right) + \mathcal{L}^f\left(\frac{1}{\log \log (n + 10)}\right)\right]\cr
\ar\ge\ar   \frac{ \Phi_\mu(n) - \Phi_\Lambda(n)}{n\log n(\log\log n)^2} + O(1)
\eeqnn
as $n \rightarrow \infty,$
which, by \eqref{Phi_0518a}, goes to $\infty.$ Then
there exists a function $d$ with $\lim_{a \rightarrow \infty}d(a) = \infty$ such that \eqref{L0208} holds for any $u \ge a$. The result follows from Proposition \ref{c0407}.
\qed

{\bf Proof of Theorem \ref{c0527a}.}
By Proposition \ref{p0217}, as $n \rightarrow \infty$ we have
\beqnn
\mathcal{L}^c\log n
\ar=\ar\sum_{k=2}^n{n\choose k}\lambda_{n,k}[\log(n-k+1)-\log n]\cr
\ar=\ar\sum_{k=2}^n{n\choose k}\lambda_{n,k}\log\left(1-\frac{k-1}{n}\right)\cr
\ar=\ar -\frac{\Phi_\Lambda(n)}{n} + O(1).
\eeqnn
On the other hand, by Condition \ref{C1} and Propositions \ref{p0211}, \ref{mu_tail}, one sees that
\beqnn
\mathcal{L}^f\log n \ar=\ar n\sum_{k=1}^\infty\mu(k)\left[\log(n+k)-\log n\right] \cr
\ar=\ar n\sum_{k=1}^\infty\mu(k)\log\left(1+\frac{k}{n}\right)\cr
\ar=\ar n\sum_{k=1}^n\mu(k)\log\left(1+\frac{k}{n}\right) + n\sum_{k=n+1}^\infty\mu(k)\log\left(1+\frac{k}{n}\right)\cr
\ar =\ar \frac{\Phi_\mu(n)}{n} + \left(2\log 2 - \frac{1}{2}\right)b(\log n)^\alpha + o((\log n)^\alpha)
\eeqnn
as $n \rightarrow \infty.$
It follows that
\beqnn
\mathcal{L}\log n\ar=\ar\mathcal{L}^c\log n+\mathcal{L}^f\log n\cr
\ar\le\ar	
-\frac{\Phi_\Lambda(n) - \Phi_\mu(n)}{n} +\left(2\log 2 - \frac{1}{2}\right) b(\log n)^\alpha + o((\log n)^{\alpha})\cr
\ar=\ar (\log n)^{\alpha} \left[-\frac{\Phi_\Lambda(n) - \Phi_\mu(n)}{n (\log n)^\alpha}  + \left(2\log 2 - \frac{1}{2}\right)b  + o(1)\right]
\eeqnn
as $n \rightarrow \infty$.
By \eqref{Phi_0527}, 	there exists a constant $\delta > 0$ such that
\beqlb\label{eq0210}
\frac{\Phi_\Lambda(n) - \Phi_\mu(n)}{n(\log n)^{\alpha}} \ge \left(2\log 2 - \frac{1}{2}\right) b + \delta
\eeqlb
for large $n$.
Then
\beqnn
\mathcal{L}\log n \le 0
\eeqnn
for large $n$.  It follows that there exists a constant $C > 0$ such that \eqref{F-L} holds for all $n \ge 1$.
The result follows by Proposition \ref{p0607}.
\qed

{\bf Proof of Theorem \ref{c0527}.}
Recall that $\alpha > 1$. By \eqref{eq0208a} and \eqref{eq0210}, for any $m \ge 1$ one sees that
\beqnn
\mathcal{L}^m\left(1 + \frac{1}{\log\log(n + l)}\right) \ar\ge\ar \frac{(\log n)^{\alpha - 1}}{(\log\log n)^2}\left[ \frac{\Phi_\Lambda(n) - \Phi_\mu(n)}{n(\log n)^\alpha} - \left(2\log 2 - \frac{1}{2}\right)b + o(1)\right]\cr
\ar\ge\ar (\delta + o(1)) \frac{(\log n)^{\alpha - 1}}{(\log\log n)^2}
\eeqnn
as $n \rightarrow \infty,$
which goes to infinity and then \eqref{L0208} holds. Moreover, the process $N_m$ does not explode by Lemma \ref{Fl3.8}. By Proposition \ref{CDI} and letting $l \rightarrow \infty$, we have
\beqnn
\lim_{a \rightarrow \infty}\lim_{n \rightarrow \infty}\mbf{P}_n\{\tau_{a,m}^- < t\} = 1
\eeqnn
for any $t > 0$ and $m \in \mathbb{N}$. Then the result follows from Proposition \ref{l0406}.
\qed

{\bf Proof of Theorem \ref{p5}.}
Recall that $\alpha > 1.$ By Lemmas \ref{L^c4} and \ref{L^f4}, as $n \rightarrow \infty,$ one obtains that
\beqnn
\mathcal{L} \left(\frac{1}{\log\log (n+10)}\right) \ar=\ar \mathcal{L}^c\left(\frac{1}{\log\log (n+10)}\right)+\mathcal{L}^f \left(\frac{1}{\log\log (n+10)}\right)\cr
\ar\le\ar \frac{\Phi_\Lambda(n) - \Phi_\mu(n)}{(n+10)\log (n+10) (\log\log (n+10))^2}\cr
\ar\ar - \left(2\log 2 - \frac{1}{2}\right)\frac{b(\log n)^{\alpha - 1}}{(\log\log n)^2} + o\left(\frac{(\log n)^{\alpha - 1}}{(\log\log n)^2}\right)\cr
\ar=\ar \frac{\Phi_\Lambda(n) - \Phi_\mu(n)}{n\log n (\log\log n)^2}- \left(2\log 2 - \frac{1}{2}\right)\frac{b(\log n)^{\alpha - 1}}{(\log\log n)^2} \cr
\ar\ar + o\left(\frac{(\log n)^{\alpha - 1}}{(\log\log n)^2}\right) + O(1),
\eeqnn
where $O(1)$ comes by changing $n + 10$ to $n$.
Notice that $2\log 2-1/2 > 0.$ Moreover, by \eqref{Phi_0518b}, there exists a constant $\delta > 0$ such that
\beqnn
\frac{\Phi_\Lambda(n) - \Phi_\mu(n)}{n(\log n)^{\alpha}} \le \left(2\log 2 - \frac{1}{2}\right)b - \delta
\eeqnn
for large $n$. Then
\beqlb\label{eq0419}
\mathcal{L} \left(\frac{1}{\log\log (n+10)}\right)
\ar\le\ar  \frac{(\log n)^{\alpha - 1}}{ (\log\log n)^2}\left[\frac{\Phi_\Lambda(n) - \Phi_\mu(n)}{n(\log n)^{\alpha}} - \left(2\log 2 - \frac{1}{2}\right)b + o(1)\right]\cr
\ar\le\ar (-\delta + o(1)) \frac{(\log n)^{\alpha - 1}}{ (\log\log n)^2}
\eeqlb
as $n \rightarrow \infty$, which goes to $-\infty$.
Then \eqref{eq_0412} holds for large $n$.
The process stays infinite by Proposition \ref{stay_infinite1}. The result follows.
\qed

{\bf Proof of Theorem \ref{p0206a}.}
Recall that $\alpha > 1$ and \eqref{Phi_0518b} holds. Then \eqref{eq_0412} holds for large $n$ by the proof of Theorem \ref{p5}. Moreover, by \eqref{eq0419}, one sees that
\beqnn
\mathcal{L} \left(1- \frac{1}{\log\log (n+10)}\right) \ar=\ar - \mathcal{L} \left( \frac{1}{\log\log (n+10)}\right) \cr
\ar\ge\ar (\delta + o(1))  \frac{(\log n)^{\alpha - 1}}{(\log\log n)^2} \rightarrow \infty
\eeqnn
as $n \rightarrow \infty$.
Then there exists a function $d$ with $\lim_{a \rightarrow \infty}d(a) = \infty$ such that \eqref{L0208} holds for any $u \ge a$. The desired result follows by Proposition \ref{c0407}.
\qed

\vskip 1cm
{\bf Acknowledgement} The authors thank Clement Foucart for very helpful discussions. The first author thanks Concordia University where part of the paper was completed during her visit.


\begin{thebibliography}{00}


\bibitem[AN72]{AN72} K. Athreya and P. Ney (1972). {\em Branching processes.} Springer, Berlin.

\bibitem[B04]{B04} J. Berestycki (2004). Exchangeable fragmentation-coalescence processes and their equilibrium
measures. {\it Electron. J. Probab.} {\bf 9}, 770--824.

\bibitem[BBL10]{BBL10} J. Berestycki, N. Berestycki and V. Limic (2010). The $\Lambda$-coalescent speed of coming down from infinity. {\it Ann. Probab.} {\bf 38}(1), 207--233.

\bibitem[B09]{B09} N. Berestycki (2009). {\it Recent Progress in Coalescent Theory. Ensaios Matem\'{a}ticos} [Mathematical Surveys]
{\bf 16}. Sociedade Brasileira de Matem\'{a}tica, Rio de Janeiro.

\bibitem[B01]{B01} J. Bertoin (2001). Homogeneous fragmentation processes. {\it Probab. Theory Related Fields} {\bf 121}(3), 301--318.

\bibitem[B02]{B02} J. Bertoin (2002). Self-similar fragmentations. {\it Ann. Inst. H. Poincar\'{e} Probab. Statist.} {\bf 38}(3), 319--340.

\bibitem[B03]{B03} J. Bertoin (2003). The asymptotic behaviour of fragmentation processes. {\it J. Euro. Math. Soc.} {\bf 5}, 395--416.

\bibitem[B06]{B06} J. Bertoin (2006). {\it Random fragmentation and coagulation processes}. Cambridge Studies in Advanced Mathematics, {\bf 102}, Cambridge University Press, Cambridge.
	
\bibitem[C86a]{C86a} M.-F. Chen (1986). Coupling of jump processes. {\it Acta Math. Sinica, New Series} {\bf 2}, 121--136.

\bibitem[C86b]{C86b} M.-F. Chen (1986). {\it Jump Processes and Interacting Particle Systems} (in Chinese).  Beijing Normal University Press.

\bibitem[C04]{C04} M.-F. Chen (2004). {\it From Markov Chains to Non-equilibrium Particle Systems}, 2nd edn. World Scientific, River Edge, NJ.
	
\bibitem[CK11]{CK11} P. L. Chou and R. Khasminskii (2011). The method of Lyapunov function for the analysis of the absorption and explosion of Markov chains. {\it Probl. Inf. Transm.} {\bf 47}(3), 232--250.


\bibitem[F22]{F22} C. Foucart (2022).
A phase transition in the coming down from infinity of simple exchangeable fragmentation-coagulation processes. {\it Ann. Appl. Probab.} {\bf 32}(1), 632--664.

\bibitem[FZ22]{FZ22} C. Foucart and X. Zhou (2022).  On the explosion of the number of fragments in simple exchangeable fragmentation-coagulation processes.
{\it Ann. Inst. Henri Poincar\'{e} Probab. Stat.} {\bf 58}(2), 1182--1207.

\bibitem[FZ23]{FZ23} C. Foucart and X. Zhou (2023). On the boundary classification of $\Lambda$-Wright-Fisher processes with frequency-dependent selection. {\it Ann. H. Lebesgue} {\bf 6}, 493--539.

\bibitem[KPRS17]{KPRS17} A. Kyprianou, S. Pagett, T. Rogers and J. Schweinsberg (2017). A phase transition in excursions from infinity of the fast fragmentation-coalescence process. {\it Ann. Probab.} {\bf 45}, 3829--3849.

\bibitem[L05]{L05} A. Lambert (2005). The branching process with logistic growth. {\it Ann. Appl. Probab.} {\bf 15}, 1506--1535.

\bibitem[LYZ19]{LYZ19} P. Li, X. Yang and X. Zhou  (2019). A general continuous-state nonlinear branching process. {\it Ann. Appl. Prob.} {\bf 29}, 2523--2555.

\bibitem[LT15]{LT15} V. Limic and A. Talarczyk (2015). Second-order asymptotics for the block counting process in a class of regularly varying $\Lambda$-coalescents. {\it Ann. Probab.} {\bf 43}(3), 1419--1455.

\bibitem[MYZ21]{MYZ21}
S. Ma., X. Yang and X. Zhou (2021). Boundary behaviors for a class of continuous-state nonlinear branching processes in critical cases. {\it Electron. Commun. Probab.} {\bf 26}, 1--10.

\bibitem[MZ23]{MZ23} R. Ma and X. Zhou (2023). Explosion of continuous-state branching processes with competition in a L\'{e}vy environment. {\it J. Appl. Prob.} {\bf 61}, 68--81.


\bibitem[MT93]{MT93} S. P. Meyn and R. L. Tweedie (1993). Stability of Markovian processes, III: Foster-Lyapunov criteria for continuous-time processes. {\it Adv. Appl. Prob.} {\bf 25}, 518--548.

\bibitem[P99]{P99} J. Pitman (1999). Coalescents with multiple collisions. {\it Ann. Probab.} {\bf 27}(4), 1870--1902.

\bibitem[P06]{P06} J. Pitman (2006). {\it Combinatorial stochastic processes}. Lecture Notes in Mathematics, {\bf  1875}, Springer-Verlag, Berlin, Lectures from the 32nd Summer School on Probability Theory held in Saint-Flour,
July 7–24, 2002, With a foreword by Jean Picard.

\bibitem[RXYZ22]{RXYZ22} Y. Ren, J. Xiong, X. Yang and X. Zhou (2022). On the extinction-extinguishing dichotomy for a stochastic Lotka-Volterra type population dynamical system. {\it Stoch. Process. Appl.} {\bf 150}, 50--90.

\bibitem[S99]{S99} S. Sagitov (1999). The general coalescent with asynchronous mergers of ancestral lines. {\it J. Appl. Probab.} {\bf 36}(4), 1116--1125.

\bibitem[S00]{S00} J. Schweinsberg (2000). A necessary and sufficient condition for the $\Lambda$-coalescent to come down from infinity. {\it Electron. Commun. Probab.} {\bf 5}, 1--11.
\end{thebibliography}
\end{document}